\documentclass[a4paper,11pt]{amsart}
\input{Macros}
\usepackage[frenchb, english]{babel}

\newcommand{\OO}{\mathcal{O}}
\newcommand{\LL}{\mathcal{L}}
\newcommand{\NN}{\mathcal{N}}
\newcommand{\MM}{\mathcal{M}}
\newcommand{\mm}{\mathfrak{m}}

\title{Positivity of the CM line bundle for K-stable log Fanos}
\author{Quentin Posva}
\date{}

\address{\'{E}cole Polytechnique F\'{e}d\'{e}rale de Lausanne, SB MATH CAG, MA C3 595 (Bâtiment MA), Station 8, CH-1015 Lausanne, Switzerland}
\email{quentin.posva@epfl.ch}

\begin{document}

\maketitle

\begin{quote}
\textsc{Abstract.} We prove the bigness of the Chow--Mumford line bundle associated to a $\bQ$-Gorenstein family of log Fano varieties of maximal variation with uniformly K-stable general geometric fibers. This result generalizes a theorem of Codogni and Patakfalvi to the logarithmic setting.
\end{quote}

\tableofcontents

\section{Introduction}
Throughout this article, we work over an algebraically closed field $k$ of characteristic zero. 

\bigskip
The notion of K-stability originates from complex analytic geometry. It was first formulated by Tian in \cite{Tian_KE_metrics_with_positive_scalar_curvature} to study the existence of Kähler--Einstein metrics on Fano manifolds, and later expressed in algebraic terms by Donaldson in \cite{Donaldson_Scalar_curvature_and_stability_of_toric_varieties}. The connection between K-stability and birational geometry, and in particular with the Minimal Model program (MMP), was first noticed several years later: Odaka showed in \cite{Odaka_GIT_stability_of_polarized_varieties} that K-stable Fano varieties are log terminal, and Li and Xu used methods from the MMP to approach questions related to K-stability \cite{Li_Xu_Special_test_configuration_and_K-stability_of_Fano_varieties}. These were the first steps of a purely algebraic K-stability theory of polarized varieties, with a particular emphasis on the study of K-stability of Fano varieties. Equivalent definitions of K-stability were afterwards formulated in terms of Ding invariant in \cite{Berman_Kps_of_Q_Fano_varieties_admitting_KE_metrics} and of valuation theory \cite{Fujita_Valuative_criterion_for_uKs_of_Fano_varieties, Fujita_Odaka_On_the_K_stability_of_Fano_varieties}. This established a solid ground to study K-stable Fano varieties with methods of birational geometry.

The algebraic geometers' interest for K-stable Fano varieties comes, amongst other reasons, from the possibility of constructing well-behaved moduli spaces. Indeed, after the work of several authors, K-stability appeared to be an adequate global stability condition to obtain compact coarse moduli spaces of Fano varieties. To wit, several compact moduli spaces of del Pezzo K-stable surfaces were constructed \cite{Odaka_Spotti_Sun_Compact_moduli_spaces_of_dP_surfaces_and_KE_metrics}, as well as the moduli space of smoothable K-polystable Fano varieties \cite{Li_Wang_Xu_Quasi_projectivity}. These constructions however rely heavily on techniques from analytic geometry, and it was desirable to find purely algebraic constructions. This program became reality by combining the progress in the algebraic theory of K-stability mentioned above, with the recent breakthroughs in birational geometry (e.g. \cite{Birkar_Cascini_Hacon_McKernan_Existence_of_minimal_models}, \cite{Hacon_McKernan_Xu_ACC_for_lct} and \cite{Birkar_Anti_pluricanonical_systems_on_Fanos}) and in abstract moduli theory (\cite{Alper_Good_Moduli_Spaces_for_Artin_stacks} and \cite{Alper_HL_Heinloth_Existence_of_moduli_spaces_for_algebraic_stacks}). Thanks to several recent works \cite{Jiang_Boundedness, Blum_Xu_Uniqueness_of_K_polystable_degenerations, Alper_Blum_HL_Xu_Reductivity, Blum_Liu_Xu_Openess_of_K_semistability, Xu_Minimizing_valuation_is_quasi_monomial, Codogni_Patakfalvi_Positivity_of_the_CM_line_bundle_for_families_of_K-stable_klt_Fanos,Xu_Zhuang_On_positivity_of_the_CM_line_bundle,
Liu_Xu_Zhuang_Finite_generation_and_K_stability}, we have now a good understanding of the algebraic moduli functor of K-stable Fano varieties. This article contributes to the study of its compactness properties.

\bigskip
We shall now explain in more details what is known about the moduli functor of K-stable Fano varieties, and what is our contribution. We refer to \autoref{section:Preliminaries} for the relevant definitions regarding K-stability and birational geometry.


We consider the moduli functor $\mathcal{M}^{\text{Kss}}_{n,v,c}$, where $c\in \bQ_+$, sending a $k$-scheme $S$ to the set
		\begin{equation*}
		\mathcal{M}^{\text{Kss}}_{n,v,c}(S)=\begin{Bmatrix}
		\text{Families }(X,c\Delta)\to S\text{ where } (X,\Delta)\to S \text{ is a}\\
		\text{family of log pairs, and for every }t\in T \text{ the log}\\
		\text{fiber }(X_t,c\Delta_t)\text{ is a K-semistable log Fano pair}\\
		\text{of dimension }n\text{ and volume }v.
		\end{Bmatrix}
		\end{equation*}
It was conjectured that this functor 
is represented by an Artin stack of finite type over $k$ and admits a projective good moduli space $M_{n,v,c}^\text{Kps}$ (in the sense of \cite{Alper_Good_Moduli_Spaces_for_Artin_stacks}), whose closed points are in bijection with $n$-dimensional K-polystable $\bQ$-Fano varieties of volume $v$. As hinted above, this conjecture is now verified, thanks to the work of several authors: 

\begin{theorem}[\cite{Jiang_Boundedness, Blum_Xu_Uniqueness_of_K_polystable_degenerations, Alper_Blum_HL_Xu_Reductivity, Blum_Liu_Xu_Openess_of_K_semistability, Xu_Minimizing_valuation_is_quasi_monomial,
Xu_Zhuang_On_positivity_of_the_CM_line_bundle,
Liu_Xu_Zhuang_Finite_generation_and_K_stability}]\label{thm:Moduli_thm}
The moduli functor $\sM_{n,v,c}^\text{Kss}$ is an Artin stack of finite type over $k$ and admits a projective good moduli space $M_{n,v,c}^\text{Kps}$ whose $k$-points parametrize K-polystable $\bQ$-Fano varieties of dimension $n$ and volume $v$.		
\end{theorem}	

At the time the first version of the present article was written, the above theorem was not yet proved entirely: the missing parts were the properness and projectivity of $M_{n,v,c}^{\text{Kps}}$, which were latter settled through the proof of Finite Generation Conjecture in \cite{Liu_Xu_Zhuang_Finite_generation_and_K_stability}. It was also conjectured, and has been verified in full generality in \textit{op.cit.}, that the polarization on the moduli space is given by the so-called \emph{Chow--Mumford (CM) line bundle}. Our article was part of the effort, together with \cite{Codogni_Patakfalvi_Positivity_of_the_CM_line_bundle_for_families_of_K-stable_klt_Fanos} and \cite{Xu_Zhuang_On_positivity_of_the_CM_line_bundle}, to show that the CM line bundle is indeed a good candidate.

Before stating our result, let us define the CM line bundle (see also \autoref{subsection:CM line bundle}). We consider $f\colon (X,D)\to T$ a flat family of log pairs of relative dimension $n$, such that $X$ and $T$ are projective and normal, and $-(K_{X/T}+\Delta)$ is $f$-ample. We let
		$$\lambda_{f,D}:=-f_*((-(K_{X/T}+D))^{n+1}),$$
where $f_*$ is the cycle-pushforward. Then $\lambda_{f,D}$ is a $\bQ$-Cartier divisor on $T$, called the CM line bundle of the family $f\colon (X,D)\to T$. It has a good functorial behaviour (see \autoref{thm:base_change_of_CM_line_bundle}) and therefore defines a $\bQ$-line bundle $\lambda$ on $\mathcal{M}^{\text{Kss}}_{n,v,c}$. Better still, it descends to the good moduli space $M^{\text{Kps}}_{n,v,c}$ in the sense that there exists a $\bQ$-line bundle $L$ on $M_{n,v,c}^\text{Kps}$ whose pullback to $\MM_{n,v,c}^\text{Kss}$ is $\lambda$ \cite[Lemma 10.2]{Codogni_Patakfalvi_Positivity_of_the_CM_line_bundle_for_families_of_K-stable_klt_Fanos}. 

Our main result, which is a step towards the ampleness of $\lambda$, reads as follows:

\begin{theorem}\label{thm:main_thm}
Let $f\colon (X,D)\to T$ be a flat morphism of relative dimension $n$ with connected fibers from a normal projective pair to a normal projective variety, such that $-(K_{X/T}+D)$ is $\bQ$-Cartier and $f$-ample. Assume that $D$ does not contain any fibers.
	\begin{enumerate}
		\item \textsc{Bigness}: If each fiber $(X_t,D_t)$ is klt, the general geometric fibers $(X_{\bar{t}},D_{\bar{t}})$ are uniformly K-stable, and the variation of $f$ is maximal, then $\lambda_{f,D}$ is big.
		\item \textsc{Ampleness}: If all the geometric fibers $(X_{\bar{t}},D_{\bar{t}})$ are uniformly K-stable, and the variation of $f$ is maximal, then $\lambda_{f,D}$ is ample.
	\end{enumerate}
	(Here a general geometric fiber denotes the fiber along a geometric point $\Spec\overline{\Omega}\to U\subseteq X$, where $U\subseteq X$ is a dense open subset and $\overline{\Omega}$ some algebraically closed field.)
\end{theorem}

The case $D=0$ of \autoref{thm:main_thm} was proved  previously in \cite[Theorem 1.9]{Codogni_Patakfalvi_Positivity_of_the_CM_line_bundle_for_families_of_K-stable_klt_Fanos}. However that proof does not generalize to the case $D\neq 0$. The difficulty lies in that there exist non-isomorphic log Fano pairs whose underlying varieties are isomorphic, so a family of log Fano pairs $(X,D)\to T$ can be of maximal variation while the underlying family $X\to T$ is not. Thus special attention to the geometry of the boundary $D$ is required. Our strategy of proof of \autoref{thm:main_thm} is explained in \autoref{section:Overview}: it relies on a perturbative argument on the boundary.

After the first version of this article was put on ArXiv, new positivity results for the CM line bundle were proved in \cite{Xu_Zhuang_On_positivity_of_the_CM_line_bundle}. The authors introduce the notion of \emph{reduced uniform K-stability}, which generalise that of uniform K-stability, and they proved the analogue of \autoref{thm:main_thm} for families of reduced uniform K-stable log Fano pairs, see \cite[\S 7]{Xu_Zhuang_On_positivity_of_the_CM_line_bundle}. Their strategy to deal with the case $D\neq 0$ was inspired by ours.

\begin{remark}
In \autoref{thm:main_thm}, one of our assumptions is that each fiber $(X_t,D_t)$ of the family is a klt pair. This hypothesis is natural for applications to moduli space of K-stable Fano varieties, where the families we consider have klt fibers (see \autoref{thm:K_stability_and_delta_invariant}). However it might not be necessary, since in the case $D=0$ we only need the general log fiber to be klt \cite[Theorem 1.9.a]{Codogni_Patakfalvi_Positivity_of_the_CM_line_bundle_for_families_of_K-stable_klt_Fanos}.
\end{remark}


\subsection{Overview of the proof}\label{section:Overview}
The proof of the bigness statement is based on the following idea. Let $(X,D)\to T$ be a $\bQ$-Gorenstein family of log Fano pairs of maximal variation with uniformly K-stable general geometric fibers. By \cite[Theorem 1.8]{Codogni_Patakfalvi_Positivity_of_the_CM_line_bundle_for_families_of_K-stable_klt_Fanos}, we know that $\lambda_{f,D}$ is a pseudo-effective divisor. Assume that the components of $D$ are $\bQ$-Cartier. Then for a small perturbation $D^\epsilon$ of $D$, the perturbed family $(X,D^\epsilon)\to T$ has the same properties as the original one. Hence the perturbed CM line bundle $\lambda_{f,D^\epsilon}$ remains pseudo-effective. By understanding the variation of $\lambda_{f,D}$ into $\lambda_{f,D^\epsilon}$, we will deduce that $\lambda_{f,D}$ belongs to the interior of the pseudo-effective cone. If the components of $D$ are not $\bQ$-Cartier, we use techniques from the Minimal Model Program (MMP) to run a similar analysis.

\subsubsection{Curve base and $\bQ$-Cartier coefficients}\label{subsubsection:Overview_in_easy_case}
The variation of $\lambda_{f,D}$ is easy to analyse when the base $T$ is a smooth curve, and all the reduced components $D^i$ of $D$ are $\bQ$-Cartier. It follows from the definition of the CM line bundle that
		$$-\deg\lambda_{f,D}=(-K_{X/T}-D)^{n+1}, \quad n+1=\dim X.$$ 
Let $D^\epsilon=D-\sum_i\epsilon_i D^i$ be a perturbed boundary. As explained above, the divisor $\lambda_{f,D^\epsilon}$ is pseudo-effective for small values of $\epsilon$, which means $\deg\lambda_{f,D^\epsilon}\geq 0$. We calculate this degree as above:
		\begin{eqnarray*}
		-\deg \lambda_{f,D^\epsilon}&=&\left(-K_{X/T}-D+\sum_i\epsilon_iD^i\right)^{n+1}\\
		&=&(-K_{X/T}-D)^{n+1}+(n+1)\sum_i\epsilon_i(-K_{X/T}-D)^n\cdot D^i+O(\epsilon^2) \\
		& = & -\deg\lambda_{f,D} +(n+1)\sum_i\epsilon_i(-K_{X/T}-D)^n\cdot D^i+O(\epsilon^2)
		\end{eqnarray*}
Hence for small values of $\epsilon$, the function $\epsilon\mapsto -\deg \lambda_{f,D^\epsilon}$ can be approximated by an affine polynomial with linear coefficients $(-K_{X/T}-D)^n\cdot D^i$. Assume that these first order derivatives $(-K_{X/T}-D)^n\cdot D^i$ are all positive. Then $\deg \lambda_{f,D}$ cannot be too small in comparison to them, for otherwise $\deg\lambda_{f,D^\epsilon}<0$ for a small value of $\epsilon$. 

We estimate these first order derivatives using the so-called product trick, pioneered in the work of Viehweg \cite{Viehweg_Wk_positivity_and_additivity_of_Kodaira_dimension}. For positive integers $r_0,\dots,r_N$, let 
		$$D^{(r_\bullet)}=X^{(r_0)}\times_T \left(D^1\right)^{(r_1)}\times_T\dots\times_T \left(D^N\right)^{(r_N)}$$
and let $L$ be the Cartier divisor on $D^{(r_\bullet)}$ given by the sum of the pullbacks of $-K_{X/T}-D$ restricted to the different factors. Then the self-intersection of $L$ depends only on $(-K_{X/T}-D)^{n+1}$, $(-K_{X/T}-D)^n\cdot D^i$ and $r_i$. On the other hand, if $r_\bullet$ is suitably chosen, we can infer some positivity of $L$ from the positivity of the sheaf
		$$\det\left(f_*\OO_X(-K_{X/T}-D)\otimes \bigotimes_{i}f_*\OO_{D^i}(-K_{X/T}-D)\right).$$
	The positivity of this determinant sheaf is a consequence of the maximal variation assumption via Koll\'{a}r's ampleness lemma. From the positivity of $L$, we deduce a positive lower bound for the first-order derivatives $(-K_{X/T}-D)^n\cdot D^i$.
	
It is useful in the argument to twist $-K_{X/T}-D$ with a sufficient multiple of $f^*\lambda_{f,D}$, since we obtain a nef divisor \cite[Theorem 1.20]{Codogni_Patakfalvi_Positivity_of_the_CM_line_bundle_for_families_of_K-stable_klt_Fanos}. This replacement has technical significance, but does not affect the strategy.

\subsubsection{General case}
The CM line bundle behaves well with respect to base-change (\autoref{thm:base_change_of_CM_line_bundle}). In particular, it holds that $\lambda_{f,D}\cdot C=\deg\lambda_{f_C,D_C}$ for a smooth curve $C$ mapping to $T$. Hence if $T$ has higher dimension, we can base-change over a general curve $C$, apply the previous case and obtain $\lambda_{f,D}\cdot C>0$. However, this does not suffice to prove that $\lambda_{f,D}$ is big, as the boundary of the cone of movable curves of $T$ need not be spanned by classes of movable irreducible curves. Nevertheless, this strategy still works if we keep a precise track of the positivity. 
	\begin{enumerate}
		\item First we need to estimate the derivatives $(-K_{X_C/C}-D_C)^n\cdot D^i_C$. We can construct $D^{(r_\bullet)}$ and $L$ as before, conclude to some positivity of $L$ and base-change to $C$. However the base-change $D^{(r_\bullet)}\times_T C$ might not be flat over a general curve $C$, which creates difficulties. Thus we construct the product from a suitable birational model of $X$ (\autoref{notation:notation_for_main_thm_products}). Then we use the ampleness lemma and the product trick to estimate the derivatives (\autoref{thm:Application_of_ampleness_lemma} and \autoref{thm:Estimation_for_all_volumes}).
		\item We can garantee that these derivatives do not simultaneously go to zero when the class $[C]$ gets closer to the boundary of the movable cone (\autoref{thm:About_A}). This is done using the theory of Knudsen--Mumford expansion, which is recalled in \autoref{section:KM_expansion}.
		\item Once we have a uniform control on the derivatives, we would like to perturb the boundary $D$. However the components $D^i$ might not be $\bQ$-Cartier. Using the techniques of the MMP, we produce a birational model $W$ of $X$ on which some components become $\bQ$-Cartier, and such that the morphism $W\to T$ has good properties (see \autoref{thm:Global_small_modifications} for the precise statement). Then we are in position to perform the perturbation argument on $W$ (\autoref{section:variation_boundary}) and conclude.
	\end{enumerate}


\subsection{Organization of the paper}
In \autoref{section:Preliminaries} we gather some notations, recall the characterization of K-stability in terms of the $\delta$-invariant and discuss base-change of divisors and the definition of the CM line bundle. The statement of the ampleness lemma is recalled in \autoref{section:Ampleness_lemma}, and we gather some facts about the Knudsen-Mumford expansion in \autoref{section:KM_expansion}. In \autoref{section:Perturbation} we show how to perturb a family of pairs so that a part of the boundary becomes $\bQ$-Cartier, while the other relevant properties of the family are preserved. \autoref{section:Proof} is devoted to the proof of the main theorem, and the technical results are gathered in the appendix.

\subsection{Acknowledgments}
I would like to thank my advisor Zsolt Patakfalvi for his interest, his support and his help. I am also grateful to Giulio Codogni, Roberto Svaldi and Maciej Zdanowicz for many helpful conversations, and to the anonymous reviewer for numerous corrections and suggestions.

\section{Preliminaries}\label{section:Preliminaries}
\subsection{Notations and conventions}
We fix once and for all an algebraically closed field $k$ of characteristic $0$. Every scheme appearing in this article is a $k$-scheme, and every morphism is a $k$-morphism.

A \emph{pair} $(X,D)$ is the data of a normal variety $X$ and an effective $\bQ$-Weil divisor $D$ such that $K_X+D$ is $\bQ$-Cartier. We refer to \cite[\S 2.3]{Kollar_Mori_Birational_geometry_of_algebraic_varieties} for the definition of \emph{klt} and \emph{lc} pairs.

A pair $(X,D)$ is \emph{Fano} if $X$ is projective and $-K_X-D$ is ample. A pair $(X,D)$ is \emph{weak log Fano} if it is a klt projective pair such that $-K_X-D$ is big and nef. A pair $(X,D)$ is \emph{log Fano} if it is a klt Fano pair. We say that $X$ is \emph{$\bQ$-Fano} if $(X,0)$ is log Fano.

If $X$ is a Noetherian variety, an open subset $U\subset X$ is called \emph{big} if $X\setminus U$ has codimension at least $2$ in $X$. More generally, if $f\colon X\to T$ is a morphism, then $U$ is \emph{relatively big over $T$} if $U_t\subset X_t$ is big for every $t\in T$.

A birational proper morphism $\pi\colon Y\to X$ between projective varieties is called \emph{small} if the exceptional locus of $\pi$ has codimension at least $2$.

\begin{definition}[General movable curves]
Let $X$ be a projective variety. A \emph{smooth curve} $C\to X$ is a non-constant morphism (not necessarily an embedding) from a projective smooth curve $C$ to $X$. We say that a smooth curve $C\to X$ is a \emph{general movable curve} if it is the normalization of a general curve in a family of curves covering $X$. 
\end{definition}

Let $Z\subsetneq X$ be a proper closed subset. When fixing a general movable curve, we can always assume that it is not contained in $Z$. By \cite[11.4.C]{Lazarsfeld_Positivity_II}, a $\bQ$-Cartier divisor $D$ on $X$ is big (resp. pseudo-effective) if and only if $D\cdot C>0$ (resp. $D\cdot C\geq 0$) for every general movable curve $C\to X$.

\begin{definition}[Families of log Fano pairs]\label{notation:Gorenstein_families}
A \emph{$\bQ$-Gorenstein family of log Fano pairs} $f\colon (X,D)\to T$ is the data of a flat projective morphism $f\colon X\to T$ between normal projective varieties, and of an effective Weil $\bQ$-divisor $D$, such that
	\begin{enumerate}
		\item the fibers of $f$ are irreducible and normal,		
		\item the support of $D$ does not contain any fiber, 
		\item $(X_t,D_t)$ is klt for each $t\in T$ (the definition of the restricted divisor $D_t$ is given in \autoref{Section:Base_change}), and
		\item $-K_{X/T}-D$ is an $f$-ample $\bQ$-Cartier divisor.
	\end{enumerate}
\end{definition}

\begin{definition}[Maximal variation]\label{definition:max_variation}
Let $f\colon (X,D)\to T$ be a $\bQ$-Gorenstein family of log Fano pairs. Then $f$ has \emph{maximal variation} if there is a non-empty dense open subset $V\subset T$ such that for every point $t\in V$, the set $\{t'\in V\mid (X_t,D_t)\cong (X_{t'},D_{t'})\}$ is finite.
\end{definition}

\begin{notation}[Coefficient parts]\label{notation:divisors}
Let $X$ be a normal variety and $D$ a Weil $\bQ$-divisor on $X$. For $c\in\bQ$, the part of coefficient $c$ of $D$ is defined to be
		$$D^{=c}:=\sum_{\text{coeff}_ED=c}E$$
where the sum runs through the set of prime Weil divisors $E$ of $X$. We have $D=\sum_{c\in\bQ} cD^{=c}$. For simplicity, if $\{c\in\bQ \mid D^{=c}\neq 0\}=\{c_1,\dots,c_m\}$, we let $D^i:=D^{=c_i}$ so that $D=\sum_{i=1}^m c_iD^i$.
We will also denote by $D^i$ the corresponding reduced closed subscheme.
\end{notation}

\begin{notation}[Volumes]
Let $D$ be a $\bQ$-Cartier divisor on a proper scheme $X$. We denote its \emph{volume} by $\Vol(D)$. We refer to \cite[\S 2.2.C]{Lazarsfeld_Positivity_I} for the definition and the properties of the volume.
\end{notation}

\begin{notation}[Intersection numbers]
Let $D$ be a $\bQ$-Cartier divisor on a proper equidimensional scheme $X$ of dimension $n$. We denote by $D^n=(D \cdots  D)$ its self-intersection. If $D$ is ample, it holds that $\Vol(D)=D^n$. 

If $X$ is also reduced and $\LL$ is a line bundle on $X$, by abuse of notation we denote by $\sL$ the associated linear equivalence class of Cartier divisors (see \cite[Corollary 1.19]{Liu_AG}). Then it makes sense to write
		$$\sL^m\cdot D^{n-m}=(\underbrace{\sL \cdots  \sL}_{m \text{ times}}\cdot \underbrace{D \cdots D}_{n-m \text{ times}}).$$
If $f\colon C\to X$ is a morphism from a smooth proper curve, we write
	$$D\cdot C=\frac{1}{r}\deg_Cf^*\sO_X(rD),$$ where $r>0$ is such that $rD$ is Cartier.
\end{notation}

\begin{notation}[$m$-fold products]\label{notation:fiber_products}
Let $f\colon X\to T$ be a morphism of proper schemes. We denote by $X^{(m)}$ the $m$-times fiber product of $X$ with itself over $T$. It comes with projection morphisms $p_i\colon X^{(m)}\to X$ for $i=1,\dots,m$ and the structural morphism $f^{(m)}\colon X^{(m)}\to T$. Given a line bundle $\LL$ on $X$, or a Cartier divisor $D$ on $X$, we write
		$$\LL^{(m)}:=\bigotimes_{i=1}^m p_i^*\LL,\quad D^{(m)}:=\sum_{i=1}^m p_i^*D.$$
\end{notation}

\subsection{K-stability}\label{section:K_stability}
In this section, we recall briefly one characterization of the $\delta$-invariant for log Fano pairs, and its relation with K-stability. We refer to \cite{Fujita_Valuative_criterion_uniform_K_stability} for the algebraic definition of K-stability in terms of test configurations.

Consider a $n$-dimensional weak log Fano pair $(X,D)$. Let $E$ be a prime divisor over $X$, and $\pi\colon Y\to X$ be a smooth birational model on which $E$ appears. We can write
		$$K_Y\equiv_\text{num}\pi^*(K_X+D)+\sum_F a_{X,D}(F)F$$
where $F$ runs through the prime divisors of $Y$. The \emph{log discrepancy} of $E$ with respect to $(X,D)$ is
		$$A_{X,D}(E):=a_{X,D}(E)+1.$$
We also define the quantity
		$$S_{X,D}(E):=\frac{1}{(-K_X-D)^n}\int_0^{+\infty}\Vol(\pi^*(-K_X-D)-xE)dx.$$

\begin{definition}[Delta invariant]
Let $(X,D)$ be a log Fano pair. The \emph{$\delta$-invariant} of $(X,D)$ is given by
		$$\delta(X,D):=\inf_E\frac{A_{X,D}(E)}{S_{X,D}(E)},$$
where $E$ runs through the prime divisors over $X$.
\end{definition}

\begin{remark}
The original definition of the delta invariant in \cite{Fujita_Odaka_On_the_K_stability_of_Fano_varieties} is formulated in terms of basis-type divisors of the anti-log-canonical linear system. However, the above one is more convenient for our purpose. The equivalence between the two definitions is proved in \cite[Theorem 4.4]{Blum_Jonsson_Thresholds} in the case $D=0$, and \cite[Theorem 4.6]{Codogni_Patakfalvi_Positivity_of_the_CM_line_bundle_for_families_of_K-stable_klt_Fanos} in the general logarithmic case.
\end{remark}

The relation between the delta invariant and K-stability of log Fano pairs is given by the following theorem. In this article, we use this characterization as the definition of uniform K-stability. See \cite[Theorems 1.1 and 2.1]{Fujita_Odaka_On_the_K_stability_of_Fano_varieties} and \cite[Theorem B]{Blum_Jonsson_Thresholds} for the proof of equivalence.
\begin{theorem}\label{thm:K_stability_and_delta_invariant}
Let $(X,D)$ be a Fano pair.
	\begin{enumerate}
		\item $(X,D)$ is K-semistable if and only if $(X,D)$ is klt and $\delta(X,D)\geq 1$.
		\item $(X,D)$ is uniformly K-stable if and only if $(X,D)$ is klt and $\delta(X,D)>1$.
	\end{enumerate}
\end{theorem}

The next result will be useful in \autoref{section:Perturbation}.

\begin{proposition}\label{thm:variation_of_delta_invariant}
Let $(X,D)$ be a weak log Fano pair, and $\Gamma$ an effective $\bQ$-Cartier divisor supported on $\text{Supp}(D)$. Assume that 
		$$\inf_E\frac{A_{X,D}(E)}{S_{X,D}(E)}>1,$$
where $E$ runs through the divisors over $X$. Then for all rational $\epsilon >0$ small enough,
		$$\inf_E\frac{A_{X,D-\epsilon\Gamma}(E)}{S_{X,D-\epsilon\Gamma}(E)}>1.$$
\end{proposition}
\begin{proof}
Replacing $\Gamma$ by a small multiple, we may assume that there is an effective $\bQ$-Cartier divisor $\Gamma'$ on $X$ such that 
$\Gamma+\Gamma'\in |-K_X-D|_{\bQ}$. By assumption, there is an $a>0$ such that $A_{X,D}(E)\geq (1+a)S_{X,D}(E)$ for all divisors $E$. Choose $a'\in (0,a)$ and define the function
		$$f(\epsilon):=(1+a')(1+\epsilon)^{n+1}\frac{(-K_X-D)^n}{(-K_X-D+\epsilon\Gamma)^n},\quad \epsilon\in\bR.$$
Since $\lim_{\epsilon\to 0}f(\epsilon)=1+a'$, we can fix $\epsilon_0=\epsilon_0(a')>0$ such that for all $\epsilon\in (0,\epsilon_0)$, we have $f(\epsilon)<1+a$. Since we assumed $\text{Supp}(\Gamma)\subset\text{Supp}(D)$, we can also arrange that $D-\epsilon_0\Gamma$ is effective.
		
Now let $E$ be a divisor over $X$, appearing on a birational model $\pi\colon Y\to X$. For any $\epsilon>0$ and any $x\in\bR_+$, observe that
		\begin{eqnarray*}
		\Vol(\pi^*(-K_X-D+\epsilon\Gamma)-xE) &\leq & \Vol(\pi^*(-K_X-D+\epsilon(\Gamma+\Gamma'))-xE)\\
		&=&\Vol((1+\epsilon)\pi^*(-K_X-D)-xE).
		\end{eqnarray*}
Integrating over $x$, we obtain
	\begin{eqnarray*}
	\int_0^{\infty}\Vol(\pi^*(-K_X-D+\epsilon\Gamma)-xE)dx&\leq &\int_0^\infty\Vol((1+\epsilon)\pi^*(-K_X-D)-xE)dx\\
	&\overset{y=(1+\epsilon)x}{=}&(1+\epsilon)^{n+1}\int_0^\infty \Vol(\pi^*(-K_X-D)-yE)dy.
	\end{eqnarray*}
Together with the definition of the functional $S_{X,\bullet}$, the above inequality implies
		$$S_{X,D-\epsilon \Gamma}(E)\leq (1+\epsilon)^{n+1}\frac{(-K_X-D)^n}{(-K_X-D+\epsilon\Gamma)^n} S_{X,D}(E).$$
Now take $\epsilon<\epsilon_0$ and assume that $A_{X,D-\epsilon\Gamma}(E)\leq (1+a')S_{X,D-\epsilon\Gamma}(E)$. Then we have
		\begin{eqnarray*}
		A_{X,D}(E)\leq A_{X,D}(E)+\epsilon\cdot\text{ord}_E(\Gamma)=A_{X,D-\epsilon\Gamma}(E)&\leq& (1+a')S_{X,D-\epsilon\Gamma}(E)\\
		& \leq & f(\epsilon)S_{X,D}(E)\\
		&<&(1+a)S_{X,D}(E),
		\end{eqnarray*}
and we have obtained a contradiction with the hypothesis. Thus $A_{X,D-\epsilon\Gamma}(E)>(1+a')S_{X,D-\epsilon\Gamma}$ for all $\epsilon<\epsilon_0$. Since $\epsilon_0$ does not depend on $E$, the proof is complete.
\end{proof}

\subsection{Base-change of divisors}\label{Section:Base_change}
\begin{notation}\label{notation:Divisors_in_families}
In this subsection, we consider a flat morphism $f\colon X\to T$ between normal projective varieties, and an effective Weil $\bQ$-divisor $D$ on $X$, such that
	\begin{enumerate}
		\item the fibers of $f$ are connected and normal,
		\item the support of $D$ does not contain any fiber.
	\end{enumerate}
\end{notation}

\begin{remark}
If $(X,D)\to T$ is as in \autoref{notation:Divisors_in_families}, then $cD$ is a Mumford divisor in the sense of \cite[Definition 1]{Kollar_Families_of_divisors} for $c>0$ divisible enough. Since $T$ is reduced, $cD$ is automatically K-flat over $T$ \cite[Definition 2 and paragraph 6]{Kollar_Families_of_divisors}. Thus $(X,D=\frac{1}{c}\cdot {cD})\to T$ is a K-flat family in the sense of Koll\'{a}r.
\end{remark}

\begin{lemma}\label{lemma:boundary_dominates_the_base}
In the situation of \autoref{notation:Divisors_in_families}, each component of $D$ dominates $T$.
\end{lemma}
\begin{proof}
We may assume that $D$ is irreducible. Since $D$ does not contain any component of any fiber, the scheme-theoretic restriction $D\cap X_t$ has dimension at most $\dim f-1$ for any $t\in T$. If $D$ does not dominate $T$, then $\dim f(D)<\dim T$. But in this case, for a general point $t\in f(D)$ we have
		$$\dim X-1=\dim D=\dim D\cap X_t+\dim f(D)<\dim f-1+\dim T=\dim X-1,$$
	which is a contradiction.
\end{proof}

\begin{definition}[Divisorial pullbacks]\label{definition:base_change_of_divisors}
Let $U\subset X$ be the smooth locus of $f$. By assumption $U$ is relatively big over $T$, thus $f(U)=T$ and $U$ is a big open subset of $X$. By \cite[Theorem 4.21]{Kollar_Families_of_varieties_of_general_type}, every Weil divisor on $X$ not containing any component of a fiber is Cartier over $U$. 

Let $u\colon S\to T$ be a morphism from a normal variety $S$. We define the \emph{divisorial pullback of $D$ along $u$} as follows. The open set $U_S=U\times_T S$ is big in $X_S=X\times_S T$. Since $D|_U$ is $\bQ$-Cartier, it pullbacks to a $\bQ$-Cartier divisor $D|_{U_S}$ on $U_S$. We let the divisorial pullback $D_{X_S}$ of $D$ along $u$ be the unique Weil $\bQ$-divisor extending the $\bQ$-Cartier divisor $D|_{U_S}$.

In particular, if $t\in T$ is a closed point, then $D_{X_t}$ is the unique Weil $\bQ$-divisor of $X_t$ extending the $\bQ$-Cartier divisor $D|_{U\cap X_t}$. For ease of notation, we will mostly write $D_t:=D_{X_t}$.

It follows from this definition that there is a $\bQ$-linear equivalence
		\begin{equation}\label{eqn:base_change_of_pairs}
		K_{X_S/S}+D_{X_S}\sim_\bQ v^*(K_{X/T}+D)
		\end{equation}
where $v\colon X_S\to X$ is the induced morphism, see \cite[\S 2.4.1]{Codogni_Patakfalvi_Positivity_of_the_CM_line_bundle_for_families_of_K-stable_klt_Fanos}.
\end{definition}

\begin{lemma}\label{thm:base_change_of_Cartier_divisors}
In the situation of \autoref{notation:Divisors_in_families}, let $S\to T$ be a morphism from a normal projective variety. If $D$ is Cartier, then the divisorial pullback of $D$ and the pullback of $D$ as Cartier divisor along $\sigma\colon X_S\to X$ agree.
\end{lemma}
\begin{proof}
If $U$ is the smooth locus of $X\to T$, then $\sigma^*D$ represents $D_{X_S}$ on $U_S$ by definition. A Cartier divisor on a normal variety is determined in codimension one, and $U_S$ is big. Thus $\sigma^*D$ represent the Weil divisor $D_{X_S}$.
\end{proof}

\begin{lemma}\label{thm:Generic_properties_for_restricted_divisors_I}
In the situation of \autoref{notation:Divisors_in_families}, there is a dense big open set $U\subseteq T$ over which all the possible unions of components of $D$ (with the reduced structure) are flat.
\end{lemma}
\begin{proof}
Let $E$ be a union of components of $D$ with the reduced structure. By generic flatness, the locus of $T$ over which $E$ is flat, say $U_E$, is dense open. Pick a codimension one point $t\in T$, and any $x\in E$ such that $f(x)=t$. The morphism $\OO_{T,t}\to \OO_{E,x}$ is flat if and only if the uniformizer $\pi$ of $\OO_{T,t}$ is sent to a non-zero-divisor of $\OO_{E,x}$. Now if $\pi$ is a zero-divisor in $\OO_{E,x}$, then the components of $X_t$ passing through $x$ are contained in $E$. But by assumption $X_t$ is irreducible and $E$ does not contain any fiber, so this cannot happen. Thus $E\to T$ is flat at $x$. Since $x$ is arbitrary, we conclude that $t\in U_E$. Therefore $U=\bigcap_E U_E$ is big. 
\end{proof}

\begin{lemma}\label{thm:Generic_properties_for_restricted_divisors_II}
In the situation of \autoref{thm:Generic_properties_for_restricted_divisors_I}, there is a dense open set $V\subseteq U$ such that:
	\begin{enumerate}
		\item for any $v\in V$ and $c\in\mathbb{Q}$, the divisorial restriction $(D^{=c})_v$ is equal to the coefficient part $(D_v)^{=c}$.
		\item for any $v\in V$ and $c\in\mathbb{Q}$, the scheme-theoretic fiber $D^{=c}\times k(v)$ is equal to the divisorial restriction $(D^{=c})_v$ with the reduced structure.
	\end{enumerate}
\end{lemma}
\begin{proof}
Base-changing if necessary, we may assume that $U=T$. Given a reduced Weil divisor $E$ not containing any fiber, we claim that the divisorial restriction $E_t$ is reduced for a general $t\in T$. In view of \autoref{definition:base_change_of_divisors}, we may assume that $E$ is Cartier. Since the claim is local on $X$, we may assume that $E$ is actually principal, say cut out by $s\in \OO(X)$. Then $\OO_X/(s)$ is reduced and flat over $T$; thus its fiber $\OO_{X}/(s)\otimes k(t)$ over a general $t\in T$ is reduced \cite[12.2.1]{EGA_IV.3}. This means exactly that the divisor $E_t$ is reduced. If $E'$ is another reduced divisor not containing any fiber, such that $E$ and $E'$ have no common components, by applying the claim to $E+E'$ we see that the divisorial restrictions $E_t$ and $E'_t$ have no common component for a general $t\in T$. Let $E$ run through the coefficient parts $D^{=c}$ of $D$ to obtain the first assertion.

If we consider $E$ as a reduced closed subscheme, the scheme-theoretic fiber $E\times k(t)$ has pure codimension one for all $t\in T$ \cite[III.9.6]{Hartshorne_Algebraic_Geometry} and is reduced for a general $t\in T$. Combining this and the first assertion, we obtain the second assertion.
\end{proof}

\begin{corollary}\label{thm:Base_change_of_divisors}
In the situation of \autoref{thm:Generic_properties_for_restricted_divisors_II}, let $C\to T$ be a smooth curve whose image intersect $V$. Let $Z:=X\times_T C$ and $D_Z=\sum_{c\in\bQ} c(D_Z)^{=c}$ be the divisorial pullback of $D$. Then $(D_Z)^{=c}$ is the divisorial pullback of $D^{=c}$ for all $c\in\bQ$. 
\end{corollary}
\begin{proof}
We have to check that the divisorial pullbacks of any two distinct coefficient parts $D^{=c}$ and $D^{=c'}$, have no component in common. Since these divisorial pullbacks are horizontal over $C$, this can be checked on a general fiber of $Z\to C$. Since $C$ meets $V$, the result follows from \autoref{thm:Generic_properties_for_restricted_divisors_II}.
\end{proof}

\subsection{The CM line bundle}\label{subsection:CM line bundle}
Let $f\colon X\to T$ be a flat projective morphism of relative dimension $n$ between normal projective varieties, let $D$ be an effective $\bQ$-divisor on $X$ such that $-(K_{X/T}+D)$ is $\bQ$-Cartier and $f$-ample. Assume also that the fibers of $f$ are irreducible and normal, and that $\text{Supp}(D)$ does not contain any fiber. Then the \emph{Chow--Mumford line bundle} of $f\colon (X,D)\to T$ is defined to be
		$$\lambda_{f,D}:=-f_*((-K_{X/T}-D)^{n+1})$$
where $f_*$ denotes the pushforward of cycles. By \cite[Proposition 3.7]{Codogni_Patakfalvi_Positivity_of_the_CM_line_bundle_for_families_of_K-stable_klt_Fanos}, $\lambda_{f,D}$ is a $\bQ$-Cartier $\bQ$-Weil divisor. It is compatible with base-change in the following sense:
\begin{proposition}\label{thm:base_change_of_CM_line_bundle}
In the above situation, let $\tau\colon S\to T$ be a morphism from a normal variety $S$. Let $f_S\colon X_S\to S$ be the induced morphism and $D_S$ be the divisorial pullback in the sense of \autoref{Section:Base_change}. Then 
		$\tau^*\lambda_{f,D}= \lambda_{f_S,D_S}$.
\end{proposition}
We refer to \cite[\S 3]{Codogni_Patakfalvi_Positivity_of_the_CM_line_bundle_for_families_of_K-stable_klt_Fanos} for the proof and more background.

\section{Ampleness lemma}\label{section:Ampleness_lemma}

The next theorem, which is an elaboration of \cite[3.9]{Kollar_Projectivity_moduli}, will be useful to establish positivity properties of line bundles. 

Let us fix our notations for the Grassmannians and the general linear groups over $k$. Given integers $w\geq q$, we let $\Gr_k(w,q)$ be the Grassmannian of $q$-dimensional quotients of a $w$-dimensional $k$-vector space. Given an integer $n$, we let $\GL_k(n)$ be the general linear group over $k$ of degree $n$.

\begin{theorem}\label{thm:Ampleness_lemma}
Let $U$ be a normal variety that can be embedded as a big open subset of a projective variety. Let $W,Q_1,\dots, Q_s$ be vector bundles on $U$ of respective ranks $w,q_1,\dots,q_s$. Assume that there exist morphisms $\phi_i\colon W\to Q_i$ for $i=1,\dots,s$ which are generically surjective. Assume also that the classifying map 
		$$U(k)\to \left( \prod_{i=1}^s \Gr_k(w,q_i)\right) \big/ \GL_k(w)$$
		is finite-to-one on a dense subset of $U$. Then for any ample Cartier divisor $B$ on $U$, there exists a positive integer $m>0$ and a non-zero morphism
				$$\emph{Sym}^{mq}\left( \bigoplus_{i=1}^{s^2w} W\right)\longrightarrow \OO_U(-B)\otimes \left( \bigotimes_{i=1}^s \det Q_i\right)^{\otimes m}$$
		where $q=\sum_{i=1}^s q_i$.
\end{theorem}
\begin{proof}
The proof is contained in \cite[\S 5]{Kovacs_Patakfalvi_Projectivity_of_the_moduli_space_of_stable_log-varieties}. More precisely, let $Q_i'$ be the image of $\phi_i\colon W\to Q_i$, with corestricted morphisms $\phi_i'\colon W\twoheadrightarrow Q_i'$. There is a big open subset $U'\subset U$ over which $Q_i'$ is locally free of rank $q_i$, for each $i$. Since the statement is about the existence of a non-zero map between two locally free sheaves, by reflexivity we may replace $U$ by $U'$ and assume that all $Q_i'$ are locally free. Now let $W'=\oplus_{i=1}^s W$, $Q'=\oplus_{i=1}^s Q_i'$ and $\phi'=\oplus_{i=1}^s \phi_i'$. As explained in \cite[Lemma 5.6]{Kovacs_Patakfalvi_Projectivity_of_the_moduli_space_of_stable_log-varieties}, $\phi'$ is surjective over $U$, and there is a dense open set of $U$ where the classifying map corresponding to $\phi'$ has finite fibers. Now we follow the proof of \cite[Theorem 5.5]{Kovacs_Patakfalvi_Projectivity_of_the_moduli_space_of_stable_log-varieties} applied to $\phi'\colon W'\to Q'$. In this proof, the assumption of weak positivity is only used in the last three lines of proof; in particular, the equation $(5.5.5)$ and the subsequent displayed isomorphisms hold without this assumption. Thus they apply to our setting: for some $m>0$, there exists a non-zero morphism
		$$\text{Sym}^{m\cdot \rk(Q')}\left(\bigoplus_{i=1}^{\rk(W')} W' \right)\longrightarrow \OO_U(-B)\otimes \left(\det Q' \right)^{\otimes m}.$$
Moreover, it follows from \autoref{lemma:det_preserves_generically_surjective_inclusions} below that we have an inclusion 
	$$\sO_U(-B)\otimes (\det Q')^{\otimes m}\hookrightarrow \sO_U(-B)\otimes \left( \bigotimes_{i=1}^s \det Q_i\right)^{\otimes m}.$$
The result follows by composing the two morphisms.
\end{proof}

\begin{lemma}\label{lemma:det_preserves_generically_surjective_inclusions}
Let $U$ be a normal variety and $\alpha\colon Q'\hookrightarrow Q$ an inclusion of locally free sheaves. Assume that $\alpha$ is generically surjective. Then $\det(Q')\hookrightarrow \det(Q)$.
\end{lemma}
\begin{proof}
Let $\sF$ be the cokernel of $\alpha$. Then $\sF$ is torsion, and the determinant sheaf $\det(\sF)=\det(Q)\otimes\det(Q')^{-1}$ is of the form $\sO_U(E)$ for an effective divisor $E$. It follows that $\det(Q')\cong \sO_U(-E)\otimes\det(Q)$ embeds into $\sO_U\otimes \det(Q)\cong \det(Q)$.
\end{proof}

\section{About the Knudsen--Mumford expansion}\label{section:KM_expansion}
In this section, we give an alternative description of the CM line bundle, that will be useful to study its positivity properties.

To begin with, we recall a special case of \cite[Theorem 4]{Knudsen_Mumford}. Let $f\colon X\to T$ be a projective morphism between Noetherian schemes of relative dimension $n$. We do not require that $f$ is flat. Let $A$ be an $f$-very ample Cartier divisor. Then there exist $\sM_i\in\text{Pic}(T)$ such that for every $q\gg 0$, there is an isomorphism
		$$\det f_*\OO_X(qA)\cong \det Rf_*\OO_X(qA)\cong \bigotimes_{i=0}^{n+1}\sM_i^{\otimes \binom{q}{i}}$$
We call this expression the \emph{Knudsen-Mumford expansion} of $\OO_X(A)$, and refer to the $\sM_{i}$ as the coefficients of the expansion. This isomorphism is moreover functorial: if $S\to T$ is a morphism from a Noetherian scheme and $A_S$ the pullback of $A$ to $X\times_TS$, then it holds that
		$$\det (f_S)_*\OO_{X_S}(qA_S)\cong \bigotimes_{i=0}^{n+1}(\sM_i)_S^{\otimes \binom{q}{i}}, \quad q\gg 0.$$
		
Now consider the particular case where $f\colon (X,D)\to T$ is a $\bQ$-Gorenstein family of log Fano pairs.
Let $s$ be such that $s(-K_{X/T}-D)$ is very ample over $T$. Then one can show that
		$$-s^{n+1}\lambda_{f,D}=\sM_{n+1},$$
where $\sM_{n+1}$ is the leading coefficient of the Knudsen-Mumford expansion of $s(-K_{X/T}-D)$. See \cite[Proposition 3.7]{Codogni_Patakfalvi_Positivity_of_the_CM_line_bundle_for_families_of_K-stable_klt_Fanos} for a proof.
	
\bigskip	
The following proposition characterizes the numerical class of $\sM_{n+1}$ in several situations: when $A$ is nef, or when the morphism $f$ has pleasant properties.

\begin{proposition}\label{thm:Leading_term_of_KM_intersection}
Let $f\colon X\to T$ be an equidimensional projective morphism of relative dimension $n$ between Noetherian proper schemes (we do not require $f$ to be flat). Let $A$ be an $f$-very ample Cartier divisor on $X$, and $\sM_{n+1}$ be the leading coefficient of the Knudsen-Mumford expansion of $\OO_X(A)$. 
	\begin{enumerate}
		\item Assume that $A$ is nef. For any smooth curve $C\to T$, it holds that
		$\sM_{n+1}\cdot C=A_C^{n+1}.$
		\item Assume that $A$ is nef and $X$ is generically reduced. Let $X'\to X$ be the normalization morphism, $f'\colon X\to T$ be the induced morphism and $A'$ be the pullback of $A$. Then $\sM_{n+1}\cdot C= f'_*((A')^{n+1})\cdot C$ for a general movable curve $C\to T$.
		\item Assume that $T$ is normal and $f$ is flat with normal fibers. Then for any smooth curve $C\to T$, we have $\sM_{n+1}\cdot C= f_*(A^{n+1})\cdot C=A_C^{n+1}$. 
	\end{enumerate} 
\end{proposition}
\begin{proof}
Fix a smooth curve $C\to T$. In any case, since both $qA$ and $qA_C$ are relatively very ample for $q\gg 0$, both sheaves $f_*\sO_X(qA)$ and $(f_C)_*\OO_{X_C}(qA_C)$ are locally free with vanishing $R^i$, $i>0$. It follows from the functoriality of the Knudsen-Mumford expansion that
		\begin{equation}\label{eqn:KM_over_curve}
		\det\left[f_*\sO_X(qA)\right]\cdot C=\deg \det\left[(f_C)_*\sO_{X_C}(qA_C)\right], \quad q\gg 0.
		\end{equation}
With this set-up:
	\begin{enumerate}
		\item Assume that $A$ is nef. The left-hand side of equation \autoref{eqn:KM_over_curve} is given by 
	$$\sum_{i=0}^{n+1}\binom{q}{i}\sM_i\cdot C=\frac{q^{n+1}}{(n+1)!}\sM_{n+1}\cdot C +O(q^n),$$ 
where $\sM_i$ are the Knudsen-Mumford coefficients of $\sO_X(A)$. Now consider the right-hand side of the same equation. By Riemann-Roch, for $q$ large enough we have
		$$h^0(X,\OO_{X_C}(qA_C))=\deg \det\left[(f_C)_*\OO_{X_C}(qA_C)\right]+\chi(C,\OO_C)\cdot \text{rk }(f_C)_*\OO_{X_C}(qA_C).$$
Since $A_C$ is nef, we have
		$$h^0(X, \sO_{X_C}(qA_C))=\frac{q^{n+1}}{(n+1)!}A_C^{n+1}+O(q^n)$$
by \cite[VI.2.15]{Kollar_Rational_curves}. Since the fibers of $f$ are $n$-dimensional, the function 
	$$q\mapsto \text{rk }(f_C)_*\OO_{X_C}(qA_C)$$ 
is a polynomial in $q$ of degree at most $n$. Hence
		$$\deg \det\left[(f_C)_*\OO_{X_C}(qA_C)\right]=\frac{q^{n+1}}{(n+1)!}A_C^{n+1}+O(q^n).$$
It follows by comparing the leading coefficients in \autoref{eqn:KM_over_curve} that $A_C^{n+1}=\sM_{n+1}\cdot C$. 
	\item Assume that $A$ is nef and $X$ generically reduced. The normalization morphism $X'\to X$ is finite, so $A'$ is nef and relatively ample over $T$. Say that $sA'$ is relatively very ample for some $s>0$, and let $\sM'_{n+1}$ be the leading coefficient of the Knudsen-Mumford polynomial of $\OO_{X'}(sA')$. Since $X$ is generically reduced, the normalization $X'\to X$ is an isomorphism away from a closed subset $Z\subsetneq X$. If $C\to T$ is a smooth curve which intersects $f(X-Z)$, the pullback morphism $(X')_C\to X_C$ is birational. Using the first assertion, we obtain that 
		$$\sM_{n+1}\cdot C=A_C^{n+1}=(A')_C^{n+1}=s^{-n-1}\sM'_{n+1}\cdot C.$$
By \cite[Lemma A.2]{Codogni_Patakfalvi_Positivity_of_the_CM_line_bundle_for_families_of_K-stable_klt_Fanos} it holds that $\sM'_{n+1}=f'_*(sA')^{n+1}$, so the second assertion follows.
\item Assume that $T$ is normal and $f$ is flat with normal fibers. Then both $X$ and $X_C$ are normal. It follows from \cite[Lemma A.2]{Codogni_Patakfalvi_Positivity_of_the_CM_line_bundle_for_families_of_K-stable_klt_Fanos} that for $q\gg 0$
		$$\det f_*\OO_X(qA)=\frac{q^{n+1}}{(n+1)!}f_*(A^{n+1})+O(q^n)$$
and
		$$\det (f_C)_*\OO_{X_C}(qA_C)=\frac{q^{n+1}}{(n+1)!}(f_C)_*(A^{n+1}_C)+O(q^n)$$
in the Chow groups of $T$ and $C$ respectively.	It follows that $\sM_{n+1}=f_*(A^{n+1})$, and by intersecting with $C$ that
		$$\sM_{n+1}\cdot C=f_*(A^{n+1})\cdot C=\deg (f_C)_*(A_C^{n+1})=A_C^{n+1}$$
		as claimed.
	\end{enumerate}

\end{proof}

\section{Perturbation of families of K-stable log Fanos}\label{section:Perturbation}
Consider a $\bQ$-Gorenstein family $(X,D)\to T$ of log Fano pairs of maximal variation, with uniformly K-stable general fibers. We show in \autoref{thm:Global_small_modifications} below that we can find a model of $(X,D)$ with the same properties over $T$, and on which some components of the boundary $D$ become $\bQ$-Cartier.

\bigskip
We need a few preliminary lemmas. For the first one, we use the terminology of \cite{Birkar_Cascini_Hacon_McKernan_Existence_of_minimal_models}.

\begin{lemma}\label{thm:stability_of_log_canonical_model}
Let $f\colon X\to T$ be a projective morphism between quasi-projective normal varieties. Let $D$ be an effective $\bQ$-Cartier $\bQ$-divisor on $X$ such that $(X,D)$ is (klt) weak log Fano over $T$. Assume that $D_1,\dots,D_m$ are effective $\bQ$-Cartier $\bQ$-divisors on $X$ with supports contained in the support of $D$. Then there exists a full-dimensional closed polytope $P\subset(\bR_{\geq 0})^m$ with the following properties:
	\begin{enumerate}
		\item $P$ contains the origin, and its interior $\emph{int}(P)$ is contained in $(\bR_{>0})^m$; and
		\item for every rational vector $(\epsilon_1,\dots,\epsilon_m)\in \emph{int}(P)$, the log canonical models of $(-K_X-D+\sum_i \epsilon_iD_i)$ over $T$ have isomorphic underlying varieties.
	\end{enumerate}
\end{lemma}
\begin{proof}
Fix a general very ample divisor $A$ on $X$, which has no component in common with $D$. Since $X$ is of Fano type over $T$, there is a $\bQ$-boundary $\Delta$ such that $(X,\Delta)$ is klt, $\Delta$ is big over $T$ and $a(-K_X-D-A)\sim_{\bQ,T} K_X+\Delta$ for some small rational $a>0$. Replacing $A$ by a general member of its linear system, we may assume that $(X,\Delta+aA)$ is also klt \cite[5.17]{Kollar_Mori_Birational_geometry_of_algebraic_varieties}. We have $K_X+\Delta+aA+\sum_i\epsilon_i D_i\sim_{\bQ,T}a(-K_X-D)+\sum_i \epsilon_iD_i$. So a log canonical model of $(X,\Delta+aA+\sum_i\epsilon_i D_i)$ over $T$ is also a $(-K_X-D+\sum_i \frac{\epsilon_i}{a}D_i)$-log canonical model of $X$ over $T$. Therefore it is equivalent to prove that: there is a full-dimensional closed polytope $P\subset (\bR_{>0})^m$ containing the origin, such that for all $(\epsilon_1,\dots,\epsilon_m)\in \text{int}(P)$, the pairs $(X,\Delta+aA+\sum_i\epsilon_i D_i)$ have a log canonical model over $T$ with isomorphic underlying varieties.

Let us write $A':=aA$ and define the affine cone $V:=\Delta+\sum_i\bR_+ D_i$ in $\text{Weil}(X)_\bR$. Since $(X,\Delta+A')$ is klt, there is an open Euclidean neighborhood $\sU$ of $\Delta\in V$ such that for all $\Gamma\in \sU$, the pair $(X,\Gamma+A')$ is klt. Also, since $K_X+\Delta+A'\sim_{\bQ,T}a(-K_X-D)$ is big over $T$, we may shrink $\sU$ so that $K_X+\Gamma+A'$ is big over $T$ for all $\Gamma\in \sU$. With the notations of \cite[1.1.4]{Birkar_Cascini_Hacon_McKernan_Existence_of_minimal_models}, this implies that $\sU\subset \sE_{A',f}(V)$. 

It follows from \cite[Corollary 1.1.5 and Theorem E]{Birkar_Cascini_Hacon_McKernan_Existence_of_minimal_models} that there are finitely birational contractions $\psi_i\colon X\dashrightarrow Z_i$ over $T$, $i=1,\dots,n$, and a decomposition
	$$\sE_{A',f}(V)=\bigcup_{i=1}^n \sW_i$$
where each $\sW_i=\sW_{\psi_i,A',f}(V)$ is a rational polytope, such that for each $\Gamma\in \sW_i$, the underlying variety of a weak log canonical model of $(X,\Gamma+A')$ over $T$ is isomorphic to $Z_i$.

By \cite[Theorem 1.2]{Birkar_Cascini_Hacon_McKernan_Existence_of_minimal_models}, for every $\Gamma\in \sU$ the pair $(X,\Gamma+A')$ has a log canonical model over $T$. Relative log canonical models are in particular relative weak log canonical models. So we obtain that for any $\Gamma\in \sU\cap \sW_i$, the underlying variety of a log canonical model of $(X,\Gamma+A')$ is isomorphic to $Z_i$. 

Since $\sE_{A',f}(V)$ contains the open neighborhood of the origin of $V$, there must be a polytope $\sW_i$ which is of full dimension in $V$ and whose closure contains the origin of $V$. Thus we may find a closed full-dimensional polytope $P\subset (\bR_{\geq 0})^m$ containing the origin, with non-empty interior $\text{int}(P)\subset(\bR_{>0})^m$ such that
		$$\left(\bigcup_{(\epsilon_1,\dots,\epsilon_m)\in \text{int}(P)} \Delta +\sum_i\epsilon_i D_i\right)\subset \sU\cap \sW_i\quad \text{for some }i.$$
		 This finishes the proof.
\end{proof}

\begin{lemma}\label{thm:Stability_of_maximal_variation}
Let $(X,D)\to T$ be a $\bQ$-Gorenstein family $(X,D)\to T$ of log Fano pairs of maximal variation. Write $D=\sum_i c_iD^i$ as in \autoref{notation:divisors}. Then there is a rational number $r>0$ such that for all $i$ such that $D^i$ is $\bQ$-Cartier, and for all rational $\epsilon\in (-r;r)$, the family $(X,D+\epsilon D^i)\to T$ has maximal variation.
\end{lemma}
\begin{proof}
Take $r=\min_{i\neq j}\{\frac{1}{2}|c_i-c_j|\}$.
\end{proof}

\begin{lemma}\label{thm:Proprieties_of_global_pair}
Let $f\colon (X,D)\to T$ be a flat equidimensional proper morphism from a normal pair to a smooth variety. Assume that every fiber $(X_t,D_t)$ is klt. Then:
	\begin{enumerate}
		\item $(X,D)$ is klt, and
		\item for any closed point $t\in T$, if $H_1,\dots,H_d$ $(d=\dim T)$ are general Cartier divisors in a base-point free linear system such that in a neighbourhood of $t$ we have $\bigcap_i H_i=\{t\}$, then $(X,D+\sum_if^*H_i)$ is dlt in a neighborhood of $X_t$.
	\end{enumerate}
\end{lemma}
\begin{proof}
Let $t\in T$ be a closed point, and let $H_1,\dots,H_d$ $(d=\dim T)$ be general Cartier divisors on $T$ such that $\bigcap_iH_i=\{t\}$.

To begin with, we prove that $(X,D)$ is klt. Indeed, we can choose the $H_1,\dots,H_{d-1}$ in a general linear system, so the iterated hyperplane sections $X^m:=\bigcap_{i=1}^mf^*H_i$ are normal varieties for $m\leq d-1$ \cite{Seidenberg_Hyperplane_sections_of_normal_varieties}. By inversion of adjunction, since $(X_t,D_t)$ is assumed to be klt, we obtain that $(X^{d-1},D|_{X^{d-1}}+f^*H_d)$ is plt along $X_t$ \cite[Theorem 5.50]{Kollar_Mori_Birational_geometry_of_algebraic_varieties}. Hence $(X^{d-1},D|_{X^{d-1}})$ is klt along $X_t$. We repeat this argument to obtain that $(X,D)$ is klt along $X_t$. The choice of $t$ was arbitrary, so we conclude that $(X,D)$ is klt.


Now we think of $t\in T$ as a fixed point and claim that if the $H_i$ are suitably chosen, then the pair $(X,D+\sum_{i=1}^df^*H_i)$ is dlt in a neighborhood of $X_t$. Indeed, we can choose the $H_m$ inductively with the property that 
				$$\text{for each } I\subseteq \{1,\dots,m\}, \text{ the intersection }X^I:=\bigcap_{i\in I}f^*H_i \; \text{ is irreducible and normal.}$$
			Each of these conditions is satisfied for a general $H_m$ passing through $t$, except for the condition on $H_d$ that $\bigcap_{i=1}^d f^*H_i=X_t$ is irreducible and normal. But this is satisfied for any choice of $H_d$, since $X_t$ is assumed to be irreducible and normal. 
			
			It follows from the choices of these $H_i$ that $(X,D+\sum_{i=1}^df^*H_i)$ is snc at every generic point of the $X^I$. If $E$ is an exceptional divisor over $X$ whose center $c_X(E)$ belongs to $X_t$ but does not belong to the snc locus of $(X,D+\sum_if^*H_i)$, then $c_X(E)$ defines a point of codimension $\geq 1$ in $X_t$, and by adjunction we obtain that $(X_t,D_t)$ is not klt, which is a contradiction.		
\end{proof}

\begin{proposition}\label{thm:Global_small_modifications}
Let $f\colon (X,D)\to T$ be a $\bQ$-Gorenstein family of log Fano pairs of maximal variation with uniformly K-stable general geometric fiber. Assume that $T$ is smooth. Then there exists a positive number $r_{X,D}>0$ with the following property. For every coefficient part $\Gamma:=D^i$ of $D$ (as in \autoref{notation:divisors}), there exists a small proper birational morphism $\nu\colon W\to X$ such that: 
	\begin{enumerate}
		\item the strict transform $\Gamma_W$ of $\Gamma$ is $\bQ$-Cartier, and 
		\item for any rational $0<\epsilon <r_{X,D}$, the family $(W,D_W-\epsilon \Gamma_W)\to T$ is a $\bQ$-Gorenstein family of log Fano pairs of maximal variation with uniformly K-stable general geometric fibers. (Here $D_W$ denotes the strict transform of $D$.)
	\end{enumerate}
\end{proposition}
\begin{proof}
Since there are finitely many coefficient part of $D$, we only need to prove the result for a fixed $\Gamma$. First we construct $\nu\colon W\to X$.
\begin{itemize}
		\item The pair $(X,D)$ is klt by \autoref{thm:Proprieties_of_global_pair}, so by \cite[Corollary 1.37]{Kollar_Singularities_of_the_minimal_model_program} there is a small proper birational morphism $\mu\colon Y\to X$ where $Y$ is a $\mathbb{Q}$-factorial projective variety. Denote by $D_Y$ the strict transform of $D$, and $\Gamma_Y$ the strict transform of $\Gamma$. We have
				$$\mu^*(K_X+D)\sim_\bQ K_Y+D_Y.$$
		\item For $\epsilon>0$, run a $(-K_Y-D_Y+\epsilon\Gamma_Y)$-MMP over $X$ to obtain a relative log canonical model. By \autoref{thm:stability_of_log_canonical_model}, this model $W$ is the same for all $0<\epsilon\ll 1$. Denote by $p\colon Y\dashrightarrow W$ the induced morphism and $D_W:=p_*D_Y,\Gamma_W:=p_*\Gamma_Y$. 
	\end{itemize}
\noindent Our construction is pictured by the following diagram:
	$$\begin{tikzcd}
		(Y,D_Y,\Gamma_Y)\arrow[d, "\mu"] \arrow[r, dashrightarrow, "p"] & (W,D_W,\Gamma_W)\arrow[dl, "\nu"]\arrow[ddl, "g", bend left] \\
		(X, D,\Gamma) \arrow[d, "f"] &&\\
		T&&
	\end{tikzcd}$$
We must show that for smal rationall $\epsilon >0$, the morphism $(W,D_W-\epsilon\Gamma_W)\to T$ is flat between normal projective varieties, of maximal variation, with (klt) log Fano fibers and uniformly K-stable general geometric fibers, and that $\Gamma_W$ is $\bQ$-Cartier. First we establish the global properties of $W$ and $\Gamma_W$.
	\begin{itemize}
		\item The morphism $\nu\colon W\to X$ is small. Indeed, $\mu$ is small and $p$ extracts no divisors. Notice that $D_W$ is equal to the strict transform of $D$.
		\item Since $W$ is the end product of an MMP, it is normal. Moreover, since $\nu$ is small and $(X,D)$ is klt by \autoref{thm:Proprieties_of_global_pair}, $(W,D_W)$ is klt and hence $W$ is Cohen-Macaulay.
		\item The $\bQ$-divisor $\Gamma_W$ is $\bQ$-Cartier. Indeed, $p_*(-K_Y-D_Y+\epsilon \Gamma_Y)=-K_W-D_W+\epsilon\Gamma_W$ is $\bQ$-Cartier by construction. It holds that $\nu^*(-K_X-D)=-K_W-D_W$, and so $-K_W-D_W$ is also $\bQ$-Cartier. Therefore $\Gamma_W$ is $\bQ$-Cartier.
		\item We have
				$$\epsilon \Gamma_Y\equiv_X -K_Y-(D_Y-\epsilon\Gamma_Y)$$
		and thus
				$$\epsilon \Gamma_W\equiv_X -K_W-(D_W-\epsilon\Gamma_W) \quad \text{is ample over }X.$$
		Hence $\Gamma_W$ is a $\mathbb{Q}$-Cartier $\mathbb{Q}$-divisor which is ample over $X$. Furthermore,
				\begin{eqnarray*}
					-K_W-(D_W-\epsilon\Gamma_W) & = & p_*(-K_Y-(D_Y-\epsilon\Gamma_Y)) \\
					& = & p_*(p^*\nu^*(-K_X-D)+\epsilon\Gamma_Y) \\
					& = & \nu^*(-K_X-D)+\epsilon\Gamma_W.
				\end{eqnarray*}
		Now by \cite[Lemma 2.4]{Patakfalvi_Semi_negativity_of_Hodge_bundles_associated_to_Du_Bois_families},
				$$-K_W+g^*K_T=-K_{W/T}, \quad -K_X+f^*K_T=-K_{X/T}.$$		
		Hence
				$$-K_{W/T}-(D_W-\epsilon\Gamma_W)=\nu^*(-K_{X/T}-D)+\epsilon\Gamma_W$$
		and therefore $-K_{W/T}-D_W+\epsilon\Gamma_W$ is ample over $T$ for $0<\epsilon\ll 1$. 
	\end{itemize}

\noindent To study the properties of a fiber $W_t$, we may shrink $T$ and work in a neighborhood of $W_t$. Let $H_1,\dots,H_d$ be general Cartier divisors such that $\bigcap_i H_i=\{t\}$. 
	\begin{itemize}
	\item  The fibers of $\nu$ are connected by construction. Since the fibers of $f$ are irreducible, we deduce that $g$ has connected fibers.
		\item Since $\nu$ is small, it gives a crepant morphism from $(W,D_W+\sum_{i=1}^dg^*H_i)$ to $(X,D+\sum_{i=1}^df^*H_i)$. Moreover $\nu$ is an isomorphism above the snc locus of $(X,D+\sum_if^*H_i)$ \cite[1.40]{Debarre_Higher_dimensional_algebraic_geometry}. By \autoref{thm:Proprieties_of_global_pair} the pair $(X,D+\sum_{i}f^*H_i)$ is dlt. Thus $(W,D_W+\sum_{i=1}^dg^*H_i)$ is lc, and every lc center is contained in the locus where $\nu$ is an isomorphism. It follows that $(W,D_W+\sum_{i=1}^dg^*H_i)$ is also dlt. 
		\item From \cite[4.16]{Kollar_Singularities_of_the_minimal_model_program} we deduce that every irreducible component of $W_t=\bigcap_{i=1}^dg^*H_i$ is normal of codimension $d$, and an lc center of $(W,D_W+\sum_{i=1}^dg^*H_i)$. Assume that $W_t$ has at least two different components; by connectedness of the fibers of $W\to T$, the two components must intersect, and the intersection is a union of lc centers of $(W,D_W+\sum_{i=1}^dg^*H_i)$ \cite[4.20.2]{Kollar_Singularities_of_the_minimal_model_program}. Since $(X,D)$ is klt by \autoref{thm:Proprieties_of_global_pair}, we have $\lfloor D\rfloor=0$ and thus $\lfloor D_W\rfloor =0$. So by \cite[4.16.1]{Kollar_Singularities_of_the_minimal_model_program} the components of $W_t$ are minimal lc centers, and we have reached a contradiction. Hence $W_t$ is irreducible and normal of codimension $d$.		
		\item Since $W$ is Cohen-Macaulay, $T$ smooth and the fibers $W_t$ equidimensional, the morphism $g\colon W\to T$ is flat \cite[Theorem 23.1]{Matsumura_Commutative_Ring_Theory}.
		\item Assume that some fiber $W_t$ is contained in the support of $D_W$. Since $W_t$ dominates $X_t$ and $D_W$ dominates $D$, we obtain that $X_t$ is contained in the support of $D$, which is impossible. Thus $\text{Supp}(D_W)$ contains no fibers of $W\to T$.
	\end{itemize}
Finally we study the pairs $(W_t,(D_W-\epsilon\Gamma_W)_t)$.
	\begin{itemize}
		\item We know that $(W,D_W+\sum_ig^*H_i)$ is dlt with reduced boundary $\sum_ig^*H_i$. Combining adjunction and \cite[4.16]{Kollar_Singularities_of_the_minimal_model_program}, we obtain that the pair $(W_t,(D_W)_t)$ is klt. Hence $(W_t,(D_W-\epsilon\Gamma_W)_t)$ is klt for every $0<\epsilon< \coeff_{\Gamma_W}D_W$ \cite[2.27]{Kollar_Mori_Birational_geometry_of_algebraic_varieties}.
		\item Since $\mu$ is small, for a general $t\in T$ the morphism $W_t\to X_t$ is small and $(D_W)_t$ is the strict transform of $D_t$. Fix one such $t$ for which $(X_t,D_t)$ is uniformly K-stable. Then for every prime divisor $E$ over $W_t$, we have
				$$A_{W_t,(D_W)_t}(E)=A_{X_t,D_t}(E),\quad S_{W_t,(D_W)_t}(E)=S_{X_t,D_t}(E),$$
	so by \autoref{thm:K_stability_and_delta_invariant} we have
		$$\inf_E\frac{A_{W_t,(D_W)_t}(E)}{S_{W_t,(D_W)_t}(E)}=\inf_E\frac{A_{X_t,D_t}(E)}{S_{X_t,D_t}(E)}=\delta(X_t,D_t)>1.$$
	Moreover $(W_t,(D_W)_t)$ is a weak log Fano pair. Thus we may apply \autoref{thm:variation_of_delta_invariant} to obtain that 
			$$\delta(W_t,(D_W)_t-\epsilon\Gamma_t)=\inf_E\frac{A_{W_t,(D_W)_t-\epsilon\Gamma_t}(E)}{S_{W_t,(D_W)_t-\epsilon\Gamma_t}(E)}>1,$$
	for all $\epsilon>0$ small enough depending on $t$. Since $(W_t,(D_W)_t-\epsilon\Gamma_t)$ is a log Fano pair, we conclude by \autoref{thm:K_stability_and_delta_invariant} that $(W_t,(D_W)_t-\epsilon\Gamma_t)$ is uniformly K-stable for all $\epsilon=\epsilon(t)>0$ small enough. 
	
	\item By openness of the uniform K-stable locus \cite[Theorem 6.8]{Blum_Liu_Openess_of_uniform_K-stability_in_families_of_Q-Fano_varieties}, we conclude that: for all rational $0<\epsilon \ll 1$, the general fiber of the family $(W,D_W-\epsilon\Gamma_W)\to T$ is uniformly K-stable.
			
		\item For a general $t\in T$, the pair $(X_t,D_t)$ is a log canonical model of $(W_t,(D_W)_t)$. Thus $(W_t,(D_W)_t)\cong (W_u,(D_W)_u)$ implies that $(X_t,D_t)\cong (X_u,D_u)$. Moreover by \autoref{thm:Stability_of_maximal_variation}, for $\epsilon$ small enough, $(W_t,(D_W-\epsilon\Gamma_W)_t)\cong (W_u,(D_W-\epsilon\Gamma_W)_u)$ if and only if $(W_t,(D_W)_t)\cong (W_u,(D_W)_u)$. Therefore $(W,D_W-\epsilon\Gamma_W)\to T$ has maximal variation for $0<\epsilon\ll 1$.
	\end{itemize}
	
\noindent This shows that $(W,D_W-\epsilon \Gamma_W)\to T$ has the required properties for all rational numbers $0<\epsilon\ll 1$.
\end{proof}

\section{Proof of the main result}\label{section:Proof}
This section is devoted to the proof of \autoref{thm:main_thm}, which we divide into several steps. See \autoref{section:Overview} for an overview of the strategy. In \autoref{section:General_notations}, we set-up the notational framework of the proof. We use the ampleness lemma in \autoref{section:Application_of_ampleness_lemma} to obtain the positivity of some relevant sheaf. The estimates of the derivatives using the product trick is obtained in \autoref{section:Estimation_derivatives}, and the perturbation argument is given in \autoref{section:variation_boundary}. 

\subsection{General notations}\label{section:General_notations}

\begin{notation}\label{notation:notation_for_main_thm_I}
Let $T$ be a smooth variety and $f\colon (X,D=\sum_{i=1}^N c_iD^i)\to T$ be a $\bQ$-Gorenstein family of log Fano pairs of maximal variation with uniformly K-stable general geometric fibers. Here $D^i$ is the part of coefficient $c_i$, see \autoref{notation:divisors}. We introduce the following additional notations, and shall use them for the rest of this section.
	\begin{enumerate}
		\item Let $n:=\dim X-\dim T$ and $v:=((-K_{X/T}-D)|_{X_t})^n$. We write $\delta:=\delta(X_{\bar{\eta}}, D_{\bar{\eta}})$ where $\eta$ is the generic point of $T$. Then $\delta$ is the value of $\delta(X_t,D_t)$ for a very general point $t\in T$ \cite[Proposition 4.15]{Codogni_Patakfalvi_Positivity_of_the_CM_line_bundle_for_families_of_K-stable_klt_Fanos}.
		\item We let $\lambda_{f,D}:=-f_*((-K_{X/T}-D)^{n+1})$ be the CM line bundle provided by the family $f\colon (X,D)\to T$.
		\item The restrictions of the morphism $f$ to the support of the $D^i$ (with the reduced structure) are denoted $f_i\colon D^i\to T$. We also write $D^0=X$ and $f_0=f$. By \autoref{lemma:boundary_dominates_the_base}, each $f_i$ is surjective. 
		\item We fix a rational number $\alpha> \max \{1,	\frac{\delta}{(\delta-1)v(n+1)}\}$. Notice that by \cite[Theorems 1.8 and 1.20]{Codogni_Patakfalvi_Positivity_of_the_CM_line_bundle_for_families_of_K-stable_klt_Fanos}, for any $\alpha'\geq \alpha$ the $\bQ$-Cartier divisor $-K_{X/T}-D+\alpha'f^*\lambda_{f,D}$ is nef.
	\end{enumerate}
\end{notation}

\begin{notation}\label{notation:notation_for_main_thm_base_change_to_curves}
In the situation of \autoref{notation:notation_for_main_thm_I}, let $\iota\colon C\to T$ be a smooth curve. Form the Cartesian square
		$$\begin{tikzcd}
		Z \arrow[r, "\sigma"]\arrow[d, "h"] & X\arrow[d, "f"] \\
		C\arrow[r, "\iota"] & T
		\end{tikzcd}$$
Note that $Z$ is normal because $h$ is flat and its fibers are normal. Let $D_Z$ be the divisorial pullback of $D$ (see \autoref{definition:base_change_of_divisors}), with coefficient parts $D^i_Z$. According to \autoref{eqn:base_change_of_pairs} and to \autoref{thm:base_change_of_CM_line_bundle}, we have
		\begin{equation}\label{eqn:base_change_for_CM_line_bundle}
		K_{Z/C}+D_Z\sim_{\bQ}\sigma^*(K_{X/T}+D), \quad \iota^*\lambda_{f,D}=\lambda_{h,D_Z}.
		\end{equation}
\end{notation}

\subsection{Application of the ampleness lemma}\label{section:Application_of_ampleness_lemma}

\begin{proposition}\label{thm:Application_of_ampleness_lemma}
In the situation of \autoref{notation:notation_for_main_thm_I}, for $q$ divisible enough, the line bundle
		$$ \bigotimes_{i\geq 0}\; \det (f_i)_*\OO_{D^i}(q(-K_{X/T}-D+2\alpha f^*\lambda_{f,D}))$$
is big on $T$.
\end{proposition}
The argument of the proof is inspired by \cite[\S 9.4]{Codogni_Patakfalvi_Positivity_of_the_CM_line_bundle_for_families_of_K-stable_klt_Fanos} and by \cite[Theorem 7.1.1]{Kovacs_Patakfalvi_Projectivity_of_the_moduli_space_of_stable_log-varieties}.
\begin{proof}
We let $V\subseteq U$ be the open subsets of $T$ given by \autoref{thm:Generic_properties_for_restricted_divisors_I} and \autoref{thm:Generic_properties_for_restricted_divisors_II}. By the maximal variation assumption, shrinking $V$ if necessary, we may assume that for any $t\in V$, there are only finitely many $t'\in V$ such that $(X_t,D_t)\cong (X_{t'},D_{t'})$.

If $r$ is a fixed integer divisible by the Cartier index of $-K_{X/T}-D+2\alpha f^*\lambda_{f,D}$, for an arbitrary $d\in\bZ$ we write 
		$$\MM_d:=\OO_X(dr(-K_{X/T}-D+2\alpha f^*\lambda_{f,D})),\quad \MM_d^{D^i}:=\MM_d|_{D^i}.$$
We choose an integer $r\geq 2$ such that for every $d>0$:
\begin{enumerate}
		\item $-dr(K_{X/T}+D)$ and $dr\alpha\cdot \lambda_{f,D}$ are Cartier; 
		\item $\MM_d$ is $f$-very ample;
		\item $R^jf_*\MM_d=0$ for all $j\geq 1$;
		\item for each $i\geq 1$: $(R^j(f_i)_*\MM_d^{D^i})|_V=0$ for all $j\geq 1$;
		\item for each $i\geq 1$: $f_*\MM_1\to (f_i)_*\MM_1^{D^i}$ is surjective on the open set $V$.
\end{enumerate}
These conditions imply that $f_*\MM_1$ and $((f_i)_*\MM_1^{D^i})|_V$ are locally free and compatible with base-change. In particular,
		\begin{enumerate}
		\addtocounter{enumi}{4}
		\item if $s:=\text{rk }f_*\MM_1$, then $s=h^0(X_t,\MM_1|_{X_t})$ for all $t\in T$.
		\end{enumerate}
We may also assume that:
		\begin{enumerate}
		\addtocounter{enumi}{5}
		\item the multiplication maps 
				$$\text{Sym}^df_*\MM_1\to f_*\MM_{d} \quad \text{and} \quad \text{Sym}^d(f_*\MM_1)|_V\to ((f_i)_*\MM_{d}^{D_i})|_V$$ 
			are surjective. 
		\end{enumerate}
Now that $r$ is chosen, we can find $d>0$ such that:
	\begin{enumerate}
	\addtocounter{enumi}{6}
		\item For all $t\in T$, the kernel $$K_t:=\ker \left[\text{Sym}^dH^0( \MM_1|_{X_t})\longrightarrow H^0(\MM_{d}|_{X_t})\right]$$ generates $\mathcal{I}_t(d)$, where $\mathcal{I}_t$ is the ideal sheaf of $X_t$ for the embedding $$\varphi_{\MM_1|_{X_t}}\colon X_t\hookrightarrow \bP^{s-1}.$$
		Here $\varphi_{\MM_1|_{X_t}}$ is only defined up to the action of $\GL_k(s)$ on the target.
		 Hence, writing $w:=\text{rk } \text{Sym}^df_*\MM_1$ and $q_0:=\text{rk } f_*\MM_{d}$, we see that the orbit of $K_t$ in $\Gr_k(w,q_0)/\GL_k(s)$ determines the projective embedding $\varphi_{\MM_1|_{X_t}}$ of $X_t$ up to linear automorphisms of $\bP^{s-1}$. 
		\item Similarly, for all $v\in V$ and $i\geq 1$, the kernel $$K_v^{D^i}:=\ker \left[\text{Sym}^dH^0(\MM_1|_{X_v})\longrightarrow H^0(\MM_{d}^{D^i}|_{(D^i)_v})\right]$$ generates $\mathcal{I}_{v,i}(d)$, where $\mathcal{I}_{v,i}$ is the ideal sheaf of $(D^i)_v$ for the embedding  
				$$\varphi_{\MM_1^{D^i}|_{(D^i)_v}}\colon (D^i)_v\hookrightarrow X_v\overset{\varphi_{\MM_1|_{X_v}}}{\xrightarrow{\hspace*{1.5cm}}} \bP^{s-1},$$
		Here $\varphi_{\MM_1^{D_i}|_{(D_i)_v}}$ is only defined up to the action of $\GL_k(s)$ on the target. 
		 Hence, writing $q_i:=\text{rk } (f_i)_*\MM_{d}^{D^i}$, we see that the orbit of $K_v^{D^i}$ in $\Gr_k(w,q_i)/\GL_k(s)$ determines the projective embedding of $(D^i)_v$ up to linear automorphisms of $\bP^{s-1}$.
	\end{enumerate}
Having choosen $r$ and $d$ with these properties, we let
		$$W:=\text{Sym}^d(f_*\MM_1)|_U, \quad Q_0:=(f_*\MM_{d})|_U, \quad Q_i :=\left((f_i)_*\MM_{d}^{D^i}\right)\Big|_U \; (i\geq 1).$$
		The sheaves $W$ and $Q_0$ are locally free over $U$ by construction. Since $D^i$ is reduced the invertible sheaf $\MM_{d}^{D^i}$ satisfies the Serre condition $S_1$, thus its sections are supported on entire components, and each component of $D^i$ dominates $T$ by \autoref{lemma:boundary_dominates_the_base}. It follows that for $i\geq 1$, the $\sO_T$-modules $(f_i)_*\MM_{d}^{D^i}$ are torsion-free. Since $T$ is normal, the $(f_i)_*\sM_d^{D^i}$ are therefore locally free at codimension one points of $T$. Since they are also locally free over $V$, we may restrict $U$, but keeping it a big open set, so that the $Q_i$ become locally free for all $i\geq 0$.

By construction there are morphisms $W\to Q_i$ for $i\geq 0$, defined over $U$, which are surjective over $V$. We claim that the corresponding classifying map
					$$\xi\colon U(k)\longrightarrow \left( \prod_{i\geq 0} \Gr_k(w,q_i)\right)\big/ \GL_k(s)\quad \text{(here }\GL_k(s)\text{ acts diagonally)}$$
			is finite-to-one over $V(k)$. Fix $t\in V(k)$. The discussion in points $(g)$ and $(h)$ above shows the following. For $t'\in V(k)$, the equality $\xi(t)=\xi(t')$ holds if and only if: $X_t$ and $X_{t'}$ have isomorphic embeddings into $\bP^{s-1}$, and under this isomorphism $(D^i)_t$ is sent to $(D^i)_{t'}$ for every $i\geq 1$. As explained at the beginning of the proof, there are only finitely many such $t'$. Thus $\xi$ is finite-to-one on $V(k)$. Hence by \autoref{thm:Ampleness_lemma}, given an ample line bundle $B$ on $T$, there is a non-zero morphism 
	\begin{equation*}
		 \text{Sym}^{mq}\left( \bigoplus_{j=1}^l \text{Sym}^d(f_*\MM_1)\right)\Big|_U \longrightarrow \OO_U(-B_U)\otimes \left( \bigotimes_{i\geq 0}\; \left(\det (f_i)_*\MM_{d}^{D^i}\right)^{\otimes m}\right)\Big|_U
		 \end{equation*}
for some integers $l,m>0$ (the precise value of $l$ is given in \autoref{thm:Ampleness_lemma}, but it is not important for our purpose). Since $U$ is a big dense open subset of $T$, and since both sides are restrictions of locally free sheaves, by reflexivity this map extends to a non-zero morphism
	\begin{equation}\label{eqn: Ampleness lemma II}
	 \text{Sym}^{mq}\left( \bigoplus_{j=1}^l \text{Sym}^d(f_*\MM_1)\right) \longrightarrow \OO_T(-B)\otimes \left( \bigotimes_{i\geq 0}\; \left(\det (f_i)_*\MM_{d}^{D^i}\right)^{\otimes m}\right)
	\end{equation}
As the right-hand side is a line bundle, this map is generically surjective.

Now let $\iota\colon C\to T$ be a general smooth curve. We use the notations of \autoref{notation:notation_for_main_thm_base_change_to_curves}. As $f_*\MM_1$ is compatible with base-change (see the beginning of the proof), we obtain
		$$\iota^*f_*\MM_1\cong h_*\OO_Z(\sigma^*rM)\cong h_*\OO_Z(r(-K_{Z/C}-D_Z+2\alpha h^*\lambda_{h,D_Z})). $$
	As $C$ is general, the general geometric fiber of $(Z,D_Z)\to C$ is uniformly K-stable. Thus by \cite[Theorem 1.20]{Codogni_Patakfalvi_Positivity_of_the_CM_line_bundle_for_families_of_K-stable_klt_Fanos}, the divisor $-K_{Z/C}-D_Z+2\alpha h^*\lambda_{h,D_Z}$ is $h$-ample and nef. Moreover, since we can write
			$$r(-K_{Z/C}-D_Z+2\alpha h^*\lambda_{h,D_Z})=K_{Z/C}+D_Z+\underbrace{(r+1)(-K_{Z/C}-D_Z+\alpha h^*\lambda_{h,D_Z})+(r-1)\alpha h^*\lambda_{h,D_Z},}_{\text{ nef and } h\text{-ample}}$$
		we may apply \cite[Proposition 6.4]{Codogni_Patakfalvi_Positivity_of_the_CM_line_bundle_for_families_of_K-stable_klt_Fanos} to obtain that $\iota^*f_*\MM_1$ is nef. Hence the pullback of $\text{Sym}^{mq}\left( \bigoplus_{j=1}^l \text{Sym}^d(f_*\MM_1)\right)$ to $C$ is also nef. By generality of $C$, the restriction of the morphism \autoref{eqn: Ampleness lemma II} to $C$ is generically surjective, so we obtain that 
			$$\OO_T(-B)\Big|_C\otimes \left( \bigotimes_{i\geq 0}\; \left(\det (f_i)_*\MM_{d}^{D^i}\right)^{\otimes m}\right)\Big|_C \text{ is nef for a general movable }C\to T.$$ 
	This shows that the line bundle $\OO_T(-B)\otimes \left( \bigotimes_{i\geq 0}\; \left(\det (f_i)_*\MM_{d}^{D^i}\right)^{\otimes m}\right)$ is pseudo-effective. Thus $\bigotimes_{i\geq 0}\det (f_i)_*\MM_{d}^{D^i}$ is big, and we conclude by letting $q=rd$.
\end{proof}

\subsection{Estimation of the derivatives}\label{section:Estimation_derivatives}

\begin{notation}\label{notation:notation_for_main_thm_II}
In the situation of \autoref{notation:notation_for_main_thm_I}:
	\begin{enumerate} 
		\item We fix a positive integer $q$ such that the divisor $q(-K_{X/T}-D+2\alpha f^*\lambda_{f,D})$ is Cartier, and write 
			$$\NN=\OO_X(q(-K_{X/T}-D+2\alpha f^*\lambda_{f,D}))\quad  \text{and}\quad \NN^{D^i}=\sN|_{D^i}.$$
		\item According to \autoref{thm:Application_of_ampleness_lemma}, we may and will choose $q$ such that 
			$$\bigotimes_{i\geq 0}\; \det (f_i)_*\sN^{D^i} \text{ is a big line bundle}.$$
	\item If $C\to T$ is a smooth curve and $Z=X\times_T C$ as in \autoref{notation:notation_for_main_thm_base_change_to_curves}, then we let $\NN_Z:=\sigma^*\NN$. Recall from \autoref{eqn:base_change_for_CM_line_bundle} that
			$\sN_Z\cong\OO_Z(q(-K_{Z/C}-D_Z+2\alpha h^*\lambda_{f,D_Z}))$, and this line bundle is nef if $C$ meets the open locus of uniformly K-stable fibers \cite[Theorem 1.20]{Codogni_Patakfalvi_Positivity_of_the_CM_line_bundle_for_families_of_K-stable_klt_Fanos}.
	\end{enumerate}
\end{notation}

We aim to give a lower bound to the intersection numbers $(\sN_Z)^{\dim D^i_Z}\cdot D^i_Z$. As explained in \autoref{section:Overview}, the idea is to construct a product $D^{(r_\bullet)}$ over $T$ and then base-change over a general curve. In view of \autoref{thm:Volume_product_of_divisors_over_a_curve}, we want the pullback of $D^{(r_\bullet)}$ to be flat over that curve. Hence it would be convenient that the restricted morphisms $D^i\to T$ are flat already. To achieve this, we pass to a birational model of $X$. Unfortunately this makes the notation quite cumbersome.

\begin{notation}\label{notation:notation_for_main_thm_products}
In the situation of \autoref{notation:notation_for_main_thm_II}. Let $r_i$ be the generic rank of $(f_i)_*\NN^{D^i}$ $(i=0,\dots ,N)$.
\begin{enumerate}
	\item We let 
		$$D^{(r_\bullet)}:= \left( D^0\right)^{(r_0)}\times_T \dots \times_T \left(D^N\right)^{(r_N)}.$$
		The projection morphism from $D^{(r_\bullet)}$ to the $i^\text{th}$ $D^j$-factor is denoted by $p^{ij}\colon D^{(r_\bullet)}\to D^j$. We denote
		$$D^{(r_\bullet)}_\text{red}:=\left(D^{(r_\bullet)}\right)_\text{red},\quad D^{(r_\bullet)}_\text{norm}:=\text{normalization of } D^{(r_\bullet)}_\text{red}.$$
		We denote by $g\colon D^{(r_\bullet)}\to T$, $g_\text{red}\colon D^{(r_\bullet)}_\text{red}\to T$  and $g_\text{norm}\colon D^{(r_\bullet)}_\text{norm}\to T$ the structural morphisms. 
		
We define the line bundles
		$$\NN^{(r_\bullet)}=\bigotimes_{i,j}\left(p^{ij}\right)^*\NN^{D^j}$$
and
		$$\sN^{(r_\bullet)}_\text{red}:=\text{pullback of }\sN^{(r_\bullet)}\text{ to }D^{(r_\bullet)}_\text{red},
		\quad \sN^{(r_\bullet)}_\text{norm}:=\text{pullback of }\sN^{(r_\bullet)}\text{ to }D^{(r_\bullet)}_\text{norm}.$$
		\item Next we fix a small $\bQ$-factorial proper model $\mu\colon Y\to X$. There exists one since $X$ has a klt structure, see \cite[1.37]{Kollar_Singularities_of_the_minimal_model_program}. Denote by $D_Y=\sum_{i=1}^N c_iD_Y^i$ the strict transform of $D$. Let $\mu^i\colon D_Y^i\to D^i$ be the induced birational morphisms, with $\mu^0=\mu$. We write 
		$$D_Y^{(r_\bullet)}:= \left(D_Y^0\right)^{(r_0)}\times_T \dots \times_T \left(D_Y^N\right)^{(r_N)}.$$
	The projection morphism from $D_Y^{(r_\bullet)}$ to the $i^\text{th}$ $D^j_Y$-factor is denoted by $p^{ij}_Y\colon D_Y^{(r_\bullet)} \to D_Y^j$. 
	
	We define the line bundles
		$$\NN_Y^{D^i}:=(\mu^i)^*\NN^{D^i}\quad
\text{and}
		\quad \NN_Y^{(r_\bullet)}=\bigotimes_{i,j}\left(p^{ij}_Y\right)^*\NN_Y^{D_Y^j}.$$
			
			\item If $\iota \colon C\to T$ is a smooth curve, we denote by $Y_C$, respectively $D^i_{Y_C},D_{Y_C}^{(r_\bullet)}, \NN_{Y_C}^{(r_\bullet)}, \mu^i_C, p^{ij}_{Y_C},$ the scheme-theoretic pullbacks along $\iota$ of $Y$, respectively $D^i_Y, D_Y^{(r_\bullet)}, \NN_Y^{(r_\bullet)},\mu^i, p^{ij}_Y$. Notice that
					$$D_{Y_C}^{(r_\bullet)}=D_Y^{(r_\bullet)}\times_T C=\left(D_{Y_C}^0\right)^{(r_0)}\times_C \dots \times_C \left(D_{Y_C}^N\right)^{(r_N)}$$
			and that the projection morphisms $D_{Y_C}^{(r_\bullet)}\to D^j_{Y_C}$ are exactly the $p^{ij}_{Y_C}$. Notice also that if $\sN^{D^j_{Y_C}}$ is the pullback of $\sN^{D^j_Y}$ along $D^j_{Y_C}\to D^j_Y$, then
					$$\sN_{Y_C}^{(r_\bullet)}\cong \bigotimes_{i,j}\left(p^{ij}_{Y_C}\right)^*\sN^{D^j_{Y_C}}.$$
					
\end{enumerate}

	\noindent The construction of parts $(a)$ and $(b)$ is summarized by the following diagram, where the arrow $D_Y^{(r_\bullet)}\to D_\text{red}^{(r_\bullet)}$ exists by \autoref{thm:product_is_reduced_and_flat}\autoref{thm:product_is_reduced_and_flat_claim_3} given below.

		\begin{equation*}
			\begin{tikzcd}
			\left(D^{j}_Y, \sN^{D^j_Y}\right)\arrow[dd, "\mu^j" left] & \left(D_Y^{(r_\bullet)},\sN_Y^{(r_\bullet)}\right)\arrow[dr]\arrow[l, "p^{ij}_Y" above]\arrow[dd, "\mu^{(r_\bullet)}" left] & & \left(D^{(r_\bullet)}_\text{norm},\sN_\text{norm}^{(r_\bullet)}\right)\arrow[dl]\\
			& &  \left(D^{(r_\bullet)}_\text{red},\sN_\text{red}^{(r_\bullet)}\right)\arrow[dl] &\\
			\left(D^j,\sN^{D^j}\right) & \left(D^{(r_\bullet)},\sN^{(r_\bullet)}\right)\arrow[l, "p^{ij}" above] & &
			\end{tikzcd}
			\end{equation*}
			
\end{notation}

\noindent Next we establish some properties of these product varieties and these product line bundles.

\begin{lemma}\label{thm:product_is_reduced_and_flat}
In the situation of \autoref{notation:notation_for_main_thm_products}:
	\begin{enumerate}
		\item\label{thm:product_is_reduced_and_flat_claim_1} $D^{(r_\bullet)}$ is equidimensional over $T$, and every component dominates $T$; moreover there is a big open set of $T$ over which $D^{(r_\bullet)}$ is flat and reduced;
		\item\label{thm:product_is_reduced_and_flat_claim_2} $D_Y^{(r_\bullet)}$ is reduced, flat and equidimensional over $T$, and every components dominates $T$;
		\item\label{thm:product_is_reduced_and_flat_claim_3} $\mu^{(r_\bullet)}\colon D^{(r_\bullet)}_Y\to D^{(r_\bullet)}$ factors through $D^{(r_\bullet)}_\text{red}$;
		\item\label{thm:product_is_reduced_and_flat_claim_4} $\mu^{(r_\bullet)}$ is an isomorphism over each generic point of $D^{(r_\bullet)}$, and every component of $D^{(r_\bullet)}_Y$ dominates a component of $D^{(r_\bullet)}$;
		\item\label{thm:product_is_reduced_and_flat_claim_5} $D_{Y_C}^{(r_\bullet)}$ is flat equidimensional over $C$, and it is reduced if $C$ is general movable.
	\end{enumerate}
\end{lemma}
\begin{proof}
Assertion \autoref{thm:product_is_reduced_and_flat_claim_1} is proved in \cite[Lemma 7.11]{Kovacs_Patakfalvi_Projectivity_of_the_moduli_space_of_stable_log-varieties}, and assertion \autoref{thm:product_is_reduced_and_flat_claim_3} will follow immediately from assertion \autoref{thm:product_is_reduced_and_flat_claim_2}.

The pair $(X,D)$ is klt by \autoref{thm:Proprieties_of_global_pair}, so $(Y,D_Y)$ is klt and hence $Y$ is Cohen-Macaulay. The divisors $D^i_Y$ are $\bQ$-Cartier because $Y$ is $\bQ$-factorial, so each $D^i_Y$ is also Cohen-Macaulay \cite[5.25]{Kollar_Mori_Birational_geometry_of_algebraic_varieties}. Hence all the morphisms $D^i_Y\to T$ are flat \cite[Theorem 23.1]{Matsumura_Commutative_Ring_Theory}. This implies that the morphism $D_Y^{(r_\bullet)}\to T$ is flat. The fibers of $D_Y^{(r_\bullet)}\to T$ have the same dimension, so the morphism is equidimensional. If one component of $D^{(r_\bullet)}_Y$ does not dominate $T$, then it belongs to the non-flat locus, which is empty. The generic fiber if reduced, so \autoref{thm:Generic_fiber_is_reduced_implies_reducedness} implies that $D_Y^{(r_\bullet)}$ is reduced. This proves the assertion \autoref{thm:product_is_reduced_and_flat_claim_2}, and assertion \autoref{thm:product_is_reduced_and_flat_claim_5} is proved similarly.

To conclude, we must prove assertion \autoref{thm:product_is_reduced_and_flat_claim_4}. Let $V\subseteq X$ be the open subset over which $\mu\colon Y\to X$ is an isomorphism. Since $X$ is normal, $V$ is big. Thus $V_i:=V\cap D^i$ is dense in $D^i$ for all $i$, and $\mu^{(r_\bullet)}$ is an isomorphism over the open set $\sV:=V^{(r_0)}\times_T V_1^{(r_1)}\times_T\dots \times_T V_N^{(r_N)}\subset D^{(r_\bullet)}$. Assume that a generic point $\eta$ of $D^{(r_\bullet)}$ does not belong to $\sV$. Then $\eta$ belongs to a product of the form $\left(D^0\right)^{(r_0)}\times_T \dots \times_T \left(D^N\right)^{(r_N)}$ with one factor $D^i$ replaced by $D^i-V_i$. Such a product has dimension strictly smaller than $D^{(r_\bullet)}$ and so $\dim \overline{\{\eta\}}<\dim D^{(r_\bullet)}$, which contradicts equidimensionality. Thus $\eta\in \sV$. Similarly, a component of $D_Y^{(r_\bullet)}$ that is contracted by $\mu^{(r_\bullet)}$ must belong to a product of the form $\left(D_Y^0\right)^{(r_0)}\times_T \dots \times_T \left(D_Y^N\right)^{(r_N)}$ with one factor $D_Y^i$ replaced by $D_Y^i-f_i^{-1}V_i$. Such a product has dimension strictly smaller than $D_Y^{(r_\bullet)}$, which contradicts the equidimensionality of $D_Y^{(r_\bullet)}$. So assertion \autoref{thm:product_is_reduced_and_flat_claim_4} follows.
\end{proof}

\begin{lemma}\label{thm:Product_line_bundle_is_pseff}
In the situation of \autoref{notation:notation_for_main_thm_products},
	\begin{enumerate}
		\item $\sN_\emph{norm}^{(r_\bullet)}$ is relatively ample over $T$ and pseudo-effective;
		\item for a general movable curve $C\to T$, the line bundle $\sN_{Y_C}^{(r_\bullet)}$ is nef.
	\end{enumerate}
\end{lemma}
\begin{proof}
Let $\sU\subset T$ be a non-empty open subset with the property that for all $t\in T$, the pair $(X_t,D_t)$ is uniformly K-stable. Let $C\to T$ be a smooth curve whose image intersects $\sU$. Denote by $\sN_Z$ the pullback of $\sN$ on $Z=X\times_T C$, and $\sD^i_C=D^i\times_C T$. Then by \cite[Theorem 1.20]{Codogni_Patakfalvi_Positivity_of_the_CM_line_bundle_for_families_of_K-stable_klt_Fanos}, $\sN_Z$ and $\sN^{\sD^i_C}=\sN_Z|_{\sD^i_C}$ are nef. Thus the product line bundle
	$$\sN_{Z}^{(r_\bullet)}=\bigotimes_{i,j}\left(p^{ij}_{C}\right)^*\sN^{D^i_{C}}\quad \text{on } D^{(r_\bullet)}\times_T C$$
	is nef. By definition $\sN^{(r_\bullet)}_{Y_C}$ is a pullback of $\sN_Z^{(r_\bullet)}$, so it is also nef. This proves the second assertion.
	
	By construction, the line bundle $\sN^{(r_\bullet)}$ is relatively ample over $T$. Let $\tilde{C}\to D^{(r_\bullet)}$ be a smooth curve that meets $g^{-1}\sU$ (recall that $g\colon D^{(r_\bullet)}\to T$ is the structural morphism). If $\tilde{C}$ is contracted by $g$, then $\sN^{(r_\bullet)}\cdot\tilde{C}>0$ by relative ampleness. Otherwise, let $C\to T$ be the normalization of $g(\tilde{C})$. Then $\tilde{C}\to D^{(r_\bullet)}$ factors through $D^{(r_\bullet)}\times_T C$, on which the pullback of $\sN^{(r_\bullet)}$ is nef. Thus $\tilde{C}\cdot \sN^{(r_\bullet)}\geq 0$. This shows that $\sN^{(r_\bullet)}$ is pseudo-effective. Since $\sN^{(r_\bullet)}_\text{norm}$ is the pullback of $\sN^{(r_\bullet)}$ through the finite morphism $D^{(r_\bullet)}_\text{norm}\to D^{(r_\bullet)}$, the first assertion follows.
\end{proof}

We are now ready to estimate the intersection numbers $(\sN_Z)^{\dim D^i_Z}\cdot D^i_Z$ (where $Z=X\times_T C$). The first part of the proof is similar to the proof of \cite[7.1.1]{Kovacs_Patakfalvi_Projectivity_of_the_moduli_space_of_stable_log-varieties}. 

\begin{proposition}\label{thm:Estimation_for_all_volumes}
In the situation of \autoref{notation:notation_for_main_thm_products}, there are an ample Cartier divisor $A$ on some irreducible component $P$ of $D_Y^{(r_\bullet)}$ and a rational number $e=e(X,D,q)>0$ with the following property: for every general movable curve $C\to T$, letting $A_C$ be the pullback of $A$ to $P\times_TC$,  it holds that
		$$e(N+1)\cdot\Vol(A_C)\leq \sum_{j=0}^N\left(\NN_Z\right)^{\dim D^j_Z}\cdot D^j_Z.$$
In particular, for every general movable curve $C\to T$,
	\begin{equation}\label{eqn:Estimation_for_all_volumes}
	\exists \; j=j(C)\geq 0 \quad : \quad \left(\NN_Z\right)^{\dim D^j_Z}\cdot D^j_Z\geq e\cdot \Vol(A_C).
	\end{equation}
\end{proposition}
\begin{proof}
By \autoref{thm:product_is_reduced_and_flat}, the morphism $g\colon D^{(r_\bullet)}\to T$ is equidimensional and every component dominates $T$. Thus $g_\text{red}\colon D_\text{red}^{(r_\bullet)}\to T$ and $g_\text{norm}\colon D_\text{norm}^{(r_\bullet)}\to T$ are also equidimensional morphisms, and any components of $D^{(r_\bullet)}_\text{red}$ or $D^{(r_\bullet)}_\text{norm}$ dominates $T$.

By \autoref{thm:product_is_reduced_and_flat} and by the proof of \autoref{thm:Application_of_ampleness_lemma}, there is a big open subset $U\subseteq T$ over which $D^{(r_\bullet)}$ is flat and reduced, and the sheaves $(f_i)_*\left(\sN^{D^i}\right)|_U$ are locally free. Let $U^0:=g^{-1}U$. Since $U^0$ is reduced, it embeds as an open subset of $D_\text{red}^{(r_\bullet)}$.  Therefore the open set $U^0$ is big in $D_\text{red}^{(r_\bullet)}$ and meets every component, and so the preimage of $U^0$ in $D_\text{norm}^{(r_\bullet)}$ is big and meets every component.
	
	$$\begin{tikzcd}
	& (D_\text{norm}^{(r_\bullet)},\NN_\text{norm}^{(r_\bullet)})\arrow[d]  \\
	& (D_\text{red}^{(r_\bullet)},\NN_\text{red}^{(r_\bullet)}) \arrow[d] \\
	U^0 \arrow[d]\arrow[r, hook] \arrow[ur, hook] & (D^{(r_\bullet)} ,\NN^{(r_\bullet)})\arrow[d, "g"] &\\
	U\arrow[r, hook] & T
	\end{tikzcd}$$
Let us write $\LL:=\bigotimes_{i\geq 0} \; \det (f_i)_*\NN^{D^i}$.  On $U$, there is an embedding
		$$ \LL|_U \hookrightarrow \bigotimes_{i\geq 0}\bigotimes_{j=1}^{r_i}\;(f_i)_*\left(\NN^{D^i}|_U\right)\cong (g|_{U^0})_*\left(\NN^{(r_\bullet)}|_{U^0}\right)$$
where the first arrow is given by the natural embedding of $\det$ into the appropriate tensor power, and the isomorphism is given by \cite[Lemma 3.6]{Kovacs_Patakfalvi_Projectivity_of_the_moduli_space_of_stable_log-varieties}. By adjunction, we obtain a morphism 
		\begin{equation}\label{eqn:Pullback}		
		(g_\text{red}|_{U^0})^*\sL|_U=(g|_{U^0})^*\sL|_U\longrightarrow \NN^{(r_\bullet)}|_{U^0}=\sN_\text{red}^{(r_\bullet)}|_{U^0}.
		\end{equation}
Since $U^0$ dominates $T$, the map \autoref{eqn:Pullback} is non-zero. We may pull back this map to the normalization $D_\text{norm}^{(r_\bullet)}$. Since the preimage $U^0$ in $D_\text{norm}^{(r_\bullet)}$ is big, by reflexivity this pullback morphism extends to a non-zero morphism
		$$(g_\text{norm})^*\LL\longrightarrow \NN_\text{norm}^{(r_\bullet)}$$
which induces a non-zero map
	\begin{equation}\label{eqn:Pullback_on_normalization}
	\NN_\text{norm}^{(r_\bullet)}\otimes (g_\text{norm})^*\LL\longrightarrow\left(\NN_\text{norm}^{(r_\bullet)}\right)^{\otimes 2}.
	\end{equation}
The line bundle $\NN_\text{norm}^{(r_\bullet)}$ is relatively ample over $T$ and pseudo-effective by \autoref{thm:Product_line_bundle_is_pseff}. Moreover $(g_\text{norm})^*\LL$ is the pullback of a big divisor by \autoref{thm:Application_of_ampleness_lemma}. Thus the left-hand side of \autoref{eqn:Pullback_on_normalization} is big on every component by \autoref{thm:pseff+pullback_of_ample_is_big}. Hence $\NN_\text{norm}^{(r_\bullet)}$ is big on at least one component. This implies that $\sN^{(r_\bullet)}_\text{red}$ is big on at least one component.

By \autoref{thm:product_is_reduced_and_flat}, $\mu^{(r_\bullet)}\colon D^{(r_\bullet)}_Y\to D^{(r_\bullet)}$ is an isomorphism in codimension zero and factors through $D_\text{red}^{(r_\bullet)}$. Hence we obtain that $\NN_Y^{(r_\bullet)}$ is big on one component $P$ of $D_Y^{(r_\bullet)}$. So we may write
		$$\left(\NN_Y^{(r_\bullet)}|_P\right)^{\otimes m}\cong\OO_P(mA+mE)$$
where $A$ is ample and $E$ is effective on $P$, for some $m>0$.

Now we fix a general movable curve $C\to T$ with the following properties. Firstly, no component of its preimage in $P$ is contained in the support of $E$. Secondly, the line bundles $\sN_Z$ and $\sN_{Y_C}^{(r_\bullet)}$ are nef (see \autoref{thm:Product_line_bundle_is_pseff}). Thirdly, the induced morphisms $\mu_C^j\colon D^j_{Y_C}\to \sD^j_Z:=D^j\times_T C$ are birational for all $j$; this is acheviable since all $\mu^j\colon D^j_Y\to D^j$ are birational. Fourthly, for all $j$ the product $\sD^j_Z$ agrees with the $j$-coefficient part $D^j_Z$ of $D_Z$ in codimension one; this is achievable by combining \autoref{thm:base_change_of_Cartier_divisors} and \autoref{thm:Base_change_of_divisors}. Finally, $D_{Y_C}^{(r_\bullet)}$ is equidimensional and reduced (\autoref{thm:product_is_reduced_and_flat}.\autoref{thm:product_is_reduced_and_flat_claim_5}).

With such a curve $C\to T$ fixed, we write
		$$\OO_{P_C}(mA_C+mE_C)\cong \left(\NN_{Y_C}^{(r_\bullet)}\Big|_{P_C}\right)^{\otimes m}\cong \left(\bigotimes_{i,j}\left(p^{ij}_{Y_C}\right)^*\NN^{D^j_{Y_C}}\right)\Big|_{P_C}^{\otimes m}.$$
		Since $A$ is ample, its pullback $A_C$ is also ample. By our choice of $C$, the divisor $E_C$ is effective. Thus
				$$\Vol(A_C)\leq \Vol(A_C+E_C)=\left( \sN^{(r_\bullet)}_{Y_C}\Big|_{P_C}\right)^{\dim P_C}.$$
	It holds by equidimensionality that $\dim P_C=\dim D_{Y_C}^{(r_\bullet)}$. Using \cite[VI.2.7.3]{Kollar_Rational_curves}, we see that
			$$\left( \NN_{Y_C}^{(r_\bullet)}\right)^{\dim D_{Y_C}^{(r_\bullet)}}=\left( \sN^{(r_\bullet)}_{Y_C}\Big|_{P_C}\right)^{\dim D_{Y_C}^{(r_\bullet)}}+\sum_{P'}\left( \sN^{(r_\bullet)}_{Y_C}\Big|_{P'}\right)^{\dim D_{Y_C}^{(r_\bullet)}}$$
		where $P'$ runs through the component of $D_{Y_C}^{(r_\bullet)}$ not contained in $P_C$. Since $\sN_{Y_C}^{(r_\bullet)}$ is nef, the sum over $P'$ is non-negative, and therefore we obtain
		\begin{equation}\label{eqn:Estimation_volume_I}
		\Vol(A_C)\leq\left( \NN_{Y_C}^{(r_\bullet)}\right)^{\dim D_{Y_C}^{(r_\bullet)}}
		\end{equation}
	On the other hand, by \autoref{thm:Volume_product_of_divisors_over_a_curve} we have:
			\begin{equation}\label{eqn:Estimation_volume_II}				
				\left( \NN_{Y_C}^{(r_\bullet)}\right)^{\dim D_{Y_C}^{(r_\bullet)}}=\sum_{i=0}^Nd_i \left(\NN^{D^j_{Y_C}}\right)^{\dim D^j_{Y_C}}\prod_{j\neq i}\left(\NN^{D^j_{Y_C}}_t\right)^{\dim D^j_{Y_C}-1}
			\end{equation}
	for some rational numbers $d_i=d_i(X,D,q)>0$ and any closed point $t\in C$. The right-hand side of \autoref{eqn:Estimation_volume_II} can be simplified: observe that
			$$\left(\NN^{D^j_{Y_C}}\right)^{\dim D^j_{Y_C}}=\left(\NN^{\sD^j_Z}\right)^{\dim \sD^j_Z}=\left(\NN^{D^j_Z}\right)^{\dim D^j_Z}=\left(\NN_Z\right)^{\dim D^j_Z}\cdot D^j_Z.$$
		Indeed, the first equality holds because $\mu^i_C$ is birational, while the second equality holds because  $\sD^j_Z$ and $D^j_Z$ are equal in codimension 1 and $\NN$ has full support (see \cite[VI.2.7.3]{Kollar_Rational_curves}).
	Similarly to the previous displayed equalities, we also have
			$$\left(\NN^{D^j_{Y_C}}_t\right)^{\dim D_{Y_C}^j-1}= 
			(\sN^{D^j_{Y_C}})^{\dim D_{Y_C}^j-1}\cdot (D_{Y_C})_t =
			(\sN_Z)^{\dim D^j_Z-1}\cdot (D_Z^j)_t
			$$
	for $t\in C$ closed. Since $\NN_Z$ is the pullback of a relatively ample line bundle over $T$, the $(\sN_Z)^{\dim D^j_Z-1}\cdot (D_Z^j)_t$ are positive. Moreover $\NN_Z$ is nef so the quantities $\left(\NN_Z\right)^{\dim D^j_Z}\cdot D^j_Z$ are non-negative. Therefore, setting 
	$a=\max_j\left\{
	(\sN_Z)^{\dim D^j_Z-1}\cdot (D_Z^j)_t
	\right\}$
	 and $b=\max_i\{d_i\}$, it follows from \autoref{eqn:Estimation_volume_I} and \autoref{eqn:Estimation_volume_II} that
			$$\Vol(A_C)\leq \left( \NN_{Y_C}^{(r_\bullet)}\right)^{\dim D_{Y_C}^{(r_\bullet)}} \leq ab\sum_{i=0}^N \left(\NN_Z\right)^{\dim D^i_Z}\cdot D^i_Z.$$
		Notice that $ab$ depends only on $(X,D)$ and $q$, so we put $e=(ab(N+1))^{-1}$ to conclude.
\end{proof}

We have now a lower bound on the derivatives $\left(\NN_Z\right)^{\dim D^i_Z}\cdot D^i_Z$ in terms of the volume of the pullback divisor $A_C$. Of course, this depends on the curve $C\to T$, but it is possible to obtain some kind of uniformity. Indeed, the next lemma shows that $\Vol(A_C)$ cannot converge to zero when $[C]$ gets close to the boundary of the movable cone.

\begin{lemma}\label{thm:About_A}
Let $A$ be the ample $\bQ$-Cartier divisor on the component $P$ of $D_Y^{(r_\bullet)}$ given by \autoref{thm:Estimation_for_all_volumes}. Then there exists a big $\bQ$-Cartier divisor $\Psi$ on $T$ such that for a general movable curve $C\to T$, we have $\Vol(A_C)=\Psi\cdot C$.
\end{lemma}
\begin{proof}
By \autoref{thm:product_is_reduced_and_flat} the scheme $P$ is reduced and the morphism $P\to T$ is equidimensional, say of relative dimension $d$. So we may apply \autoref{thm:Leading_term_of_KM_intersection} to $(P,A)\to T$. Namely, let $A'$ be the pullback of $A$ to the normalization $P'$ of $P$, and let $f'\colon P'\to T$ be the induced morphism. Then for a general smooth curve $C\to T$, we have $\Vol(A_C)=A_C^{d+1}=f_*'(A')^{d+1}\cdot C$. By \autoref{thm:Cycle_pushforward_of_ample_is_big}, $f_*'(A')^{d+1}$ is big. We take $\Psi=f_*'(A')^{d+1}$.
\end{proof}

\subsection{Variation of the boundary}\label{section:variation_boundary}
\begin{notation}\label{notation:notation_for_main_thm_III}
In this subsection, we follow \autoref{notation:notation_for_main_thm_II} and let $\Psi$ be the big $\bQ$-Cartier divisor on $T$ obtained in \autoref{thm:About_A}.
\end{notation}

Given a general smooth curve $C\to T$, the inequality \autoref{eqn:Estimation_for_all_volumes} gives a lower bound for some intersection number $(\sN_Z)^{\dim D^j_Z}\cdot D^j_Z$. As explained in \autoref{section:Overview}, we wish to derive a lower bound on $ \lambda_{f,D}\cdot C$. The case $j=0$, corresponding to $D^0=X$, is the easiest.

\begin{proposition}[Case $j=0$]\label{thm:Estimate_for_i=0}
Let $C\to T$ be a smooth curve. Assume that $\left(\NN_Z\right)^{n+1}\geq e\cdot (\Psi\cdot C)$ for some $e>0$. Then 
		$\lambda_{f,D}\cdot C\geq e_0\cdot (\Psi\cdot C)$
for some rational number $e_0=e_0(X,D,q,e)>0$.
\end{proposition}
\begin{proof}
Recall the fact that $\sN_Z\cong\OO_Z(q(-K_{Z/C}-D_Z+2\alpha h^*\lambda_{f,D_Z}))$ (see \autoref{notation:notation_for_main_thm_II}). We have
		\begin{eqnarray*}
			\frac{1}{q^{n+1}}\left(\NN_Z\right)^{n+1} & = & (-K_{Z/C}-D_Z)^{n+1} + (n+1)2\alpha(h^*\lambda_{h,D_Z} \cdot (-K_{Z/C}-D_Z)^n) \\
			& = & \deg_C h_*((-K_{Z/C}-D_Z)^{n+1})+(n+1)2\alpha \cdot \deg\lambda_{h,D_Z} \cdot ((-K_{Z/C}-D_Z)_t)^n) \\
			& = & \deg_C \lambda_{h,D_Z} \cdot \left[-1+2\alpha v(n+1)\right] \\
			&=&( \lambda_{f,D}\cdot C)\cdot  \left[\underbrace{-1+2\alpha v(n+1)}_{>0 \text{ by choice of }\alpha}\right]
		\end{eqnarray*}
We let $e_0=e(q^{n+1}[2\alpha v(n+1)-1])^{-1}$ to obtain the desired inequality.
\end{proof}

If $j>0$, we wish to relate $(\sN_Z)^{\dim D^j_Z}\cdot D^j_Z$ to the a first-order derivative of $\lambda_{f,D}\cdot C$ as the component $D^j_Z$ is perturbed. Since $D^j_Z$ might not be $\bQ$-Cartier, we introduce a birational model where it is $\bQ$-Cartier.

\begin{notation}\label{notation:notation_for_main_thm_IV}
In the situation of \autoref{notation:notation_for_main_thm_III}. By \autoref{thm:Global_small_modifications}, we may fix $r_{X,D}\in (0;1)$ with the property that for every coefficient part $D^i$ of $D$, there is a small birational proper morphism $W_i\to X$ such that for all rational numbers $\epsilon\in (0;r_{X,D})$, the family $(W_i,D_{W_i}-\epsilon D^i_{W_i})\to T$ is a $\bQ$-Gorenstein family of log Fano pairs of maximal variation with uniformly K-stable general geometric fibers.

Fix an index $j>0$ and any smooth curve $\iota\colon C\to T$. Write $\nu\colon W:=W_j\to X$ and $\Gamma:=D^j$. We let 
		$$V:=W\times_X Z\cong W\times_T C.$$
Together with the notations of \autoref{notation:notation_for_main_thm_base_change_to_curves} we obtain the diagram
		$$\begin{tikzcd}
		(V,D_V,\Gamma_V)\arrow[dd, bend right=80, "g" left]\arrow[r, "\tau"]\arrow[d,"\mu"] & (W,D_W,\Gamma_W)\arrow[d, "\nu" left]\arrow[dd, bend left=80, "f'" right] \\
		(Z,D_Z,D^j_Z)\arrow[r, "\sigma"] \arrow[d, "h"] & (X,D,\Gamma) \arrow[d, "f" left] \\
		C \arrow[r, "\iota"] & T
		\end{tikzcd}$$
where $D_Z$ is the divisorial base-change of $D$, $D_Z^j$ the $c_j$-coefficient part of $D_Z$, and where the other $\bQ$-divisors are defined as follows:
	\begin{enumerate}
		\item let $D_W$ and $\Gamma_W$ be the $\nu$-strict transforms of respectively $D$ and $\Gamma$;
		\item let $D_V$ be defined on $V$ by the equality $K_{V/C}+D_V=\tau^*(K_{W/T}+D_W)$;
		\item let $\Gamma_V$ be defined on $V$ by $\Gamma_V=\tau^*\Gamma_W$. In particular, it is $\bQ$-Cartier.
	\end{enumerate}
Define the polynomial $F(C,j)\in \bR[t]$ by
		$$F(C,j)(t):=(-K_{V/C}-D_{V}+t\Gamma_{V})^{\dim V}=(-K_{V/C}-D_{V}+t\Gamma_{V})^{n+1}.$$
We will also use the ad hoc notation
		$$\mm_2\left(\sum_{l=0}^ra_lt^l\right):=\max_{l\geq 2}\{|a_l|\} \quad \text{where}\quad\sum_la_lt^l\in\bR[t] . $$
\end{notation}

\begin{lemma}\label{thm:Properties_of_global_perturbations_III}
In the situation of \autoref{notation:notation_for_main_thm_IV}. If $C\to T$ is general movable, then for every $0<\epsilon <r_{X,D}$, the family $(V,D_V-\epsilon\Gamma_V)\to C$ is a $\bQ$-Gorenstein family of log Fanos of maximal variation with uniformly K-stable general general geometric fibers.
\end{lemma}
\begin{proof}
For all rational $0<\epsilon<r_{X,D}$, the family $(W,D_W-\epsilon\Gamma_W)\to T$ is a $\bQ$-Gorenstein family of log Fanos of maximal variation with uniformly K-stable general geometric fibers. Since the fibers are normal and $C$ is normal, $V$ is also normal. If $C$ meets the open locus where the fibers are uniformly K-stable and of maximal variation, then the statement holds.
\end{proof}

We need to relate the intersection products one can do on $V$, to the intersection products one can do on $Z$ and on $T$. This is the purpose of the next three lemmas.

\begin{lemma}\label{thm:Properties_of_global_perturbations_I}
In the situation of \autoref{notation:notation_for_main_thm_IV},
	\begin{enumerate}
	\item $\mu^*(K_{Z/C}+D_Z)\sim_\bQ K_{V/C}+D_V$.
	\end{enumerate}
Moreover, if $C\to T$ is a general movable curve, then 
	\begin{enumerate}
	\stepcounter{enumi}
		\item $\mu\colon V\to Z$ is small birational;
		\item $D_V,\Gamma_V$ are the strict transforms of $D_Z$ and $D^j_Z$ respectively; 
		\item $\Gamma_V$ is the divisorial pullback of $\Gamma_W$ and the $c_j$-coefficient part of $D_V$.
	\end{enumerate}
\end{lemma}
\begin{proof}
By \autoref{eqn:base_change_for_CM_line_bundle} we have $\sigma^*(K_{X/T}+D)\sim_\bQ K_{Z/C}+D_Z$, and since $\nu\colon W\to X$ is small we have $K_{W/T}+D_W\sim_\bQ \nu^*(K_{X/T}+D)$. By definition of $D_V$, we obtain that $K_{V/C}+D_V\sim_\bQ \mu^*(K_{Z/C}+D_Z)$. This proves part $(a)$. 

By \autoref{thm:Global_small_modifications}, the morphism $W_t\to X_t$ is small birational for a general $t\in T$. So if $C$ meets the open locus of such $t\in T$, the morphism $\mu\colon V\to Z$ is birational and small as well. In this case $D_V,\Gamma_V$ are the strict transforms of $D_Z$ and $D^j_Z$. By \autoref{thm:Base_change_of_divisors}, if $C$ is general movable then $D^j_Z$ is the divisorial pullback of $\Gamma$, and the $c_j$-coefficient part of $D_Z$. So $\Gamma_V$ is the $c_j$-coefficient part of $D_V$ and the divisorial pullback of $\Gamma_W$. This proves parts $(b),(c)$ and $(d)$.
\end{proof}

\begin{lemma}\label{thm:Properties_of_global_perturbations_II}
In the situation of \autoref{notation:notation_for_main_thm_IV}, if $ C\to T$ is general movable, we have
		$$(-K_{V/C}-D_V)^{n+1}=(-K_{Z/C}-D_Z)^{n+1} \quad \text{and}\quad (-K_{V/C}-D_V)^n\cdot \Gamma_V= (-K_{Z/C}-D_Z)^n\cdot D^j_Z.$$
\end{lemma}
\begin{proof}
We use \autoref{thm:Properties_of_global_perturbations_I}. By part $(a)$, it holds that $K_{V/C}+D_V\sim_\bQ \mu^*(K_{Z/C}+D_Z)$. If $C\to T$ is general movable, $\mu$ is birational by part $(b)$. So the first equality follows. By part $(c)$, the morphism $\mu$ restricts to a birational morphism $\mu|_{\Gamma_V}\colon \Gamma_V\to D^j_Z$, and the second equality follows.
\end{proof}

\begin{lemma}\label{thm:About_m_2}
In the situation of \autoref{notation:notation_for_main_thm_IV}, for each $j>0$ there is a non-empty finite collection of $\mathbb{Q}$-Cartier divisors $\{\Upsilon_{j,l}\}_{l=2}^{n+1}$ on $T$ such that $\mm_2(F(C,j))=\max_{l\geq 2}\{|\Upsilon_{j,l}\cdot C|\}$ for any smooth curve $C\to T$.
\end{lemma}
\begin{proof}
Fix an index $j>0$ and let $(f'\colon W\to T, D_W,\Gamma_W)$ be as in \autoref{notation:notation_for_main_thm_IV}. The fibers of $f'$ are normal, and for small positive values of $\epsilon$, the $\bQ$-Cartier divisor $-K_{W/T}-D_W+\epsilon\Gamma_W$ is relatively ample over $T$. Thus by \autoref{thm:Leading_term_of_KM_intersection}, for any smooth curve $C\to T$ we have:
		$$f'_*(-K_{W/T}-D_W+\epsilon\Gamma_W)^{n+1}\cdot C=(-K_{V/C}-D_{V}+\epsilon\Gamma_{V})^{n+1}=F(C,j)(\epsilon). $$
Write $f'_*(-K_{W/T}-D_W+\epsilon\Gamma_W)^{n+1}=\sum_{l=0}^{n+1}\epsilon^l\Upsilon_{j,l}$ in the Chow ring of $T$. By linearity of the intersection product and of $f'_*$, one can describe $\Upsilon_{j,l}$ as a multiple of the pushforwards along $f'$ of the intersection $(-K_{W/T}-D_W)^{n+1-l}\cdot \Gamma_W^l$. We obtain
		$$\mm_2(F(C,j))=\max_{l\geq 2}\{ |\Upsilon_{j,l}\cdot C|\}$$
as claimed. Notice that the family $\{\Upsilon_{j,l}\}_{l=2}^{n+1}$ is non-empty since $n\geq 1$.
\end{proof}

We are now able to treat the case $j>0$.

\begin{proposition}[Case $j>0$]\label{thm:Estimate_for_i>0}
Let $C\to T$ be a general movable smooth curve. Assume that $\left(\NN_Z\right)^n\cdot D^j_Z\geq e\cdot(\Psi\cdot C)$ for some $j>0$ and $e>0$. Then there exists a rational number $e_1=e_1(X,D,q,e)>0$ such that
		$\lambda_{f,D}\cdot C\geq e_1\cdot (\Psi\cdot C).$
\end{proposition}
\begin{proof}
By generality $C\to T$, we may and will assume that the results of \autoref{thm:Properties_of_global_perturbations_III}, \autoref{thm:Properties_of_global_perturbations_I} and \autoref{thm:Properties_of_global_perturbations_II} hold. Thus we have
		\begin{equation}\label{eqn:Value_at_0}
		F(C,j)(0)=(-K_{Z/C}-D_Z)^{n+1}=-\deg_C\lambda_{h,D_Z}
		\end{equation}
and 
		\begin{equation}\label{eqn:Derivate_estimate_I}
		F'(C,j)(0)= (n+1)(-K_{Z/C}-D_Z)^n\cdot D^j_Z.
		\end{equation}
A direct calculation gives
		\begin{equation}\label{eqn:Explicit _form_of_volume}
		\left(\NN_Z\right)^n\cdot D^j_Z=q^n(-K_{Z/C}-D_Z)^n\cdot D^j_Z+2n\alpha\cdot \deg_C\lambda_{h,D_Z}\cdot \left(\NN^{D^j}_t\right)^{n-1}
		\end{equation}
	Combining \autoref{eqn:Derivate_estimate_I}, \autoref{eqn:Explicit _form_of_volume} and the hypothesis on $\left(\NN_Z\right)^n\cdot D^j_Z$, we obtain that
			\begin{equation}\label{eqn:Derivate_estimate_II}
			F'(C,j)(0)\geq \frac{(n+1)e}{q^n}(\Psi\cdot C) - \frac{2n(n+1)\alpha}{q^n}\left(\NN_t^{D^j}\right)^{n-1}\deg_C \lambda_{h,D_Z}.
			\end{equation}		
On the other hand, for any rational $0<\epsilon <r_{X,D}$, the family $(V,D_V-\epsilon\Gamma_V)\to C$ is a $\bQ$-Gorenstein family of log Fano pairs of maximal variation with uniformly K-stable general geometric fibers. Thus 
		\begin{equation}\label{eqn:Negativeness}
		-\deg_C\lambda_{h,D_V-\epsilon\Gamma_V}=F(C,i)(\epsilon)\leq 0 \quad \forall \; \epsilon\in (0,r_{X,D})
		\end{equation}
by \cite[Theorem 1.8.a]{Codogni_Patakfalvi_Positivity_of_the_CM_line_bundle_for_families_of_K-stable_klt_Fanos}. To conclude the proof, we are going to combine \autoref{eqn:Derivate_estimate_II} and \autoref{eqn:Negativeness} to get a negativity condition on $F(C,i)(0)=\deg \lambda_{h,D_Z}$.

For convenience, let us write
			$$\beta_0=\sup\left\{\left(\NN^{D^i}_t\right)^{n-1}\mid i>0, t\in C(k)\right\}, \quad w=\Psi.C, \quad a=\frac{(n+1)e}{q^n}, \quad b=\frac{2n(n+1)\alpha}{q^n}\beta_0$$
	 (notice that, by generic flatness and Noetherianity, $\beta_0$ is finite and actually a maximum). Therefore \autoref{eqn:Derivate_estimate_II} implies that
			\begin{equation}\label{eqn:Estimate_of_F'}
			F'(C,j)(0)\geq aw - b\deg_C \lambda_{h,D_Z}
			\end{equation}
	Assume that 
			\begin{equation}\label{eqn:Degree_CM_I}
			\deg_C\lambda_{h,D_Z}\leq \frac{a}{2b}w.
			\end{equation}
	Then the estimate \autoref{eqn:Estimate_of_F'} implies
			$$F'(C,j)(0)\geq \frac{a}{2}w>0.$$
	To summarize, we know by \autoref{eqn:Negativeness} that $F(C,j)(t)$ must be negative in a neighborhood of $t=0$, and if $F(C,j)(0)$ is small we have a positive lower bound on its first derivative. This gives an upper bound on $F(C,j)(0)$. Indeed, we apply \autoref{thm:Numerical lemma} with 
		$$G=\frac{a}{2}w, \quad H=\mm_2(F(C,j)), \quad l=r_{X,D}, \quad d=n+1\geq 2,$$
	and we obtain that $F(C,j)$ takes a strictly positive value on $[0;r_{X,D}/2)$ if
			\begin{equation}\label{eqn:Estimate_derivation_I}
			F(C,j)(0)> \max\left\{-\frac{ar_{X,D}}{2}w, -\frac{a^2}{4n}\frac{w^2}{\mm_2(F(C,j))} \right\},
			\end{equation}
		where we set $\frac{1}{\mm_2(F(C,j))}=+\infty$ if $\mm_2(F(C,j))=0$. But if \autoref{eqn:Estimate_derivation_I} holds, then we get a contradiction with \autoref{eqn:Negativeness}. Thus either \autoref{eqn:Degree_CM_I} fails, or \autoref{eqn:Degree_CM_I} holds and \autoref{eqn:Estimate_derivation_I} fails. This can be synthetized as
				\begin{equation}\label{eqn:Estimate_III}
				\deg_C \lambda_{h,D_Z}=-F(C,j)(0)\geq \min \left\{\frac{a}{2b}w, \frac{ar_{X,D}}{2}w, \frac{a^2}{4n}\frac{w^2}{\mm_2(F(C,j))} \right\}.
				\end{equation}
	 To conclude, we need to modify the right-hand side of \autoref{eqn:Estimate_III} so that the only quantity that depends on $C$ is $w=(\Psi.C)$. The only problematic term is $\frac{w^2}{\mm_2(F(C,j))}$: it can be dealt with using \autoref{thm:About_m_2}, as we explain now.
		
		Let $\{\Upsilon_{r,s}\}_{r,s}$ be the collection of $\bQ$-Cartier divisors on $T$ given by \autoref{thm:About_m_2} when considering every index $r>0$. Consider the function
				$$\overline{\text{Mov}}(T)_\bR-\{\bold{0}\}\to \bR, \quad \gamma\mapsto \frac{\max_{r,s} |\Upsilon_{r,s}\cdot\gamma|}{\Psi\cdot \gamma}.$$
		This function is well-defined since $\Psi$ is big and hence defines a strictly positive functional on $\overline{\text{Mov}}(T)_\bR-\{\bold{0}\}$. It is also continuous and invariant under $\bR^*_+$-scaling of its argument. So it admits a maximum which is strictly positive, since the numerator is not zero for all movable curves. Thus there exists $\beta_1>0$ such that 
		$$\frac{\mm_2(F(C,j))}{w}=\frac{\max_{s} |\Upsilon_{j,s}\cdot C|}{\Psi\cdot C}\leq \frac{\max_{r,s} |\Upsilon_{r,s}\cdot C|}{\Psi\cdot C} <\beta_1$$ 
	for all general movable curve $C$ and $j>0$. So \autoref{eqn:Estimate_III} implies that
				$$\deg_C\lambda_{h,D_Z}=-F(C,j)(0)\geq \min \left\{\frac{a}{2b},\frac{ar_{X,D}}{2}, \frac{a^2}{4n\beta_1} \right\}\cdot w.$$
		The quantity $e_1=\min \left\{\frac{a}{2b},\frac{ar_{X,D}}{2}, \frac{a^2}{4n\beta_1} \right\}$ depends only on $X,D,q$ and $e$. Therefore the proof is complete.
\end{proof}

\subsection{Proof of \autoref{thm:main_thm}}

\begin{proof}[Proof of point (c) of \autoref{thm:main_thm}]
Let $\tau\colon T'\to T$ be a resolution of singularities. Then the induced family $f_{T'}\colon (X_{T'},D_{T'})\to T'$ is again a $\bQ$-Gorenstein family of log Fano pairs of maximal variation with general geometric uniformly K-stable fibers. The morphism $\tau$ is birational and $\tau^*\lambda_{f,D}=\lambda_{f_{T'},D_{T'}}$ by \autoref{thm:base_change_of_CM_line_bundle}. So $\lambda_{f,D}$ is big if and only if $\lambda_{f_{T'},D_{T'}}$ is big. Thus we may assume that $T$ is smooth to begin with.
Let $C\to T$ be a general movable curve. By \autoref{thm:Estimation_for_all_volumes} and \autoref{thm:About_A}, the hypothesis of either \autoref{thm:Estimate_for_i=0} or \autoref{thm:Estimate_for_i>0} is fullfilled, with a constant $e$ that depends only on $X,D$ and $q$. Thus there is a constant $c=c(X,D,q)>0$ such that $\lambda_{f,D}\cdot C\geq c\cdot (\Psi\cdot C)$. As $\Psi$ is big, the result follows.
\end{proof}

\begin{proof}[Proof of point (d) of \autoref{thm:main_thm}]
By the Nakai-Moishezon theorem it is enough to prove that for all normal varieties $V$ mapping finitely to $T$, we have $(\lambda_{f,D}|_V)^{\dim V}>0$. Let $V'\to V$ be a resolution of singularities. By \autoref{eqn:base_change_for_CM_line_bundle} and since $(\lambda_{f,D}|_V)^{\dim V}=(\lambda_{f,D}|_{V'})^{\dim V'}$, we may replace $f\colon (X,D)\to T$ by $f_{V'}\colon (X_{V'},D_{V'})\to {V'}$. By assumption all the closed fibers of $f_{V'}$ are uniformly K-stable, hence klt. So all the fibers of $f_{V'}$ are klt. Therefore we are in position to apply point (c) of \autoref{thm:main_thm}.
\end{proof}

\section{Appendix}
We gather some technical results that are used in the text.

\begin{lemma}\label{thm:Cycle_pushforward_of_ample_is_big}
Let $f\colon X\to T$ be an equidimensional proper morphism of relative dimension $n$ between projective schemes. Assume that $T$ is smooth. Let $A$ be an ample $\bQ$-Cartier divisor on $X$. Then the cycle $f_*A^{n+1}$ is $\bQ$-Cartier and big. (Here $f_*$ denotes the cycle-theoretic pushforward.)
\end{lemma}
\begin{proof}
Since $f_*$ is linear, we may replace $A$ by a multiple and assume it is very ample. If $H_1,\dots,H_{n+1}\in |A|$ are general elements, then $f_*(H_1\cap\dots\cap H_{n+1})$ is a divisor on $T$, and it is Cartier as $T$ is smooth. Since $f_*$ preserves rational equivalence, we have $f_*(H_1\cap\dots\cap H_{n+1})\in |f_*A^{n+1}|$. It follows that this linear system is base-point free and separates points. The result now follows from \cite[2.60]{Kollar_Mori_Birational_geometry_of_algebraic_varieties}.
\end{proof}

\begin{lemma}\label{thm:pseff+pullback_of_ample_is_big}
Let $f\colon X\to T$ be proper morphism between normal projective $k$-schemes. Let $A$ be a pseudo-effective relatively ample $\bQ$-Cartier divisor on $X$, and $B$ a big $\bQ$-Cartier divisor on $T$. Then $A+f^*B$ is big on every component of $X$.
\end{lemma}
\begin{proof}
We may assume that $X$ is integral. Write $B\sim_\bQ C+E$ where $C$ is ample and $E$ effective. Fix an ample divisor $H$ on $X$. Choose $\epsilon'\in \bQ_+^*$ small enough such that $\epsilon'A+f^*C$ is ample on $X$. Then choose $\epsilon\in\bQ_+^*$ small enough such that $A+\epsilon H$ is effective, and $\epsilon'A+f^*C-(1-\epsilon')\epsilon H$ is still ample. We write
		$$A+f^*B \sim_\bQ f^*E+(1-\epsilon')(A+\epsilon H) + (\epsilon'A+f^*C-(1-\epsilon')\epsilon H)$$
so $A+f^*B$ is the sum of an effective and an ample $\bQ$-divisors. By \cite[2.60]{Kollar_Mori_Birational_geometry_of_algebraic_varieties}, it is big.
\end{proof}

\begin{lemma}\label{thm:Generic_fiber_is_reduced_implies_reducedness}
Let $f\colon X\to T$ be a flat morphism between Noetherian schemes. Assume that $T$ is integral, and that the generic fiber of $f$ is reduced. Then $X$ is reduced.
\end{lemma}
\begin{proof}
Let $x$ be an associated point of $X$. By assumption, the local morphism $\OO_{T,f(x)}\to \OO_{X,x}$ is flat. If $f(x)$ is not the generic point $\eta$ of $T$, then $\OO_{T,f(x)}$ has dimension at least one, and so its maximal ideal contains a non-zero divisor. By flatness, the image of this element is also a non-zero divisor in the maximal ideal of $\OO_{X,x}$. This contradicts the fact that $x$ is an associated point, so $f(x)=\eta$. Now $X_\eta$ is reduced, so $x$ cannot be an embedded associated point. Therefore $X$ is reduced.
\end{proof}


\begin{lemma}\label{thm:Volume_simple_product_of_divisors}
Let $X_i$ be proper schemes of dimensions $n_i$ $(i=1,\dots, r)$. Set $\mathcal{X}:=X_1\times_k\dots \times_k X_r$, with projections $p_i$ onto its factors. There is a positive rational number $c=c(n_1,\dots,n_r)$ with the following property: if $L_i$ are Cartier divisors on $X_i$ and $L:=\sum_{i=1}^rp_i^*L_i$, then
		$$L^{\dim \mathcal{X}}=c\prod_{i=1}^r L_i^{n_i}.$$
\end{lemma}
\begin{proof}
By induction on $r$, it suffices to consider the case $r=2$. In this case we have
		$$L^{n_1+n_2}=\binom{n_1+n_2}{n_1}\cdot (p_1^*L_1)^{n_1}\cdot (p_2^*L_2)^{n_2}=\binom{n_1+n_2}{n_1}\cdot L_1^{n_1}\cdot L_2^{n_2},$$
as claimed. In the general case, the precise form of the constant is
		$c=\prod_{i=1}^r \binom{\sum_{k\geq i}n_k}{n_i}$.
\end{proof}

\begin{lemma}\label{thm:Volume_product_of_divisors_over_a_curve}
Let $X_i\to T$ be flat morphisms from proper schemes of dimension $1+n_i$ to a common smooth curve $(i=1,\dots,r)$. Set $\mathcal{X}:=X_1\times_T \dots \times_T X_r$ with projections $p_i$ onto its factors. Then there are positive rational numbers $d_i=d(n_1,\dots, n_r)$  with the following property: if $L_i$ are Cartier divisors on $X_i$ and $L:=\sum_{i=1}^r p_i^*L_i$, then
		$$L^{\dim\mathcal{X}}=\sum_{i=1}^r d_i L_i^{n_i+1}\prod_{j\neq i}(L_j)_t^{n_j}$$
		for any closed $t\in T$.
\end{lemma}
\begin{proof}
Notice that $\dim \mathcal{X}=1+\sum_j n_j$. Hence $L^{\dim\mathcal{X}}$ is a weighted sum of $(p_1^*L_1)^{i_1}\cdots (p_r^*L_r)^{i_r}$ with $\sum_j i_j=1+\sum_j n_j$. Such a term is zero as soon as $i_j>1+n_j$ for some $j$. On the other hand, by the pigeon-hole principle, at least one $i_j$ is greater or equal to $1+n_j$. Thus:
		\begin{eqnarray*}
		L^{1+\sum_j n_j} &=& \sum_{i=1}^r\binom{1+\sum_j n_j}{1+n_i}(p_i^*L_i)^{1+n_i}\cdot \left( \sum_{j\neq i}p_j^*L_j\right)^{\sum_{j\neq i}n_j} \\
		& = &\sum_{i=1}^r\binom{1+\sum_j n_j}{1+n_i} L_i^{1+n_i}\cdot \left( \sum_{j\neq i}p_j^*L_j\big|_{p_i^{-1}(x_i)}\right)^{\sum_{j\neq i}n_j}
		\end{eqnarray*}
where the second equality holds for any $x_i\in X_i$ by flatness of $p_i$. Notice that the fiber of $p_i \colon \mathcal{X}\to X_i$ above $x_i$ is naturally isomorphic to the fiber product $\bigtimes_{j\neq i}(X_j)_{t_i}$ taken over $\Spec k$, where $t_i$ is the image of $x_i$ through $X_i\to T$. By flatness of $X_i\to T$, the intersection number $(L_i)_{t_i}^{n_i}$ does not depend on $t_i$. Applying \autoref{thm:Volume_simple_product_of_divisors}, we get
		$$L^{\dim \mathcal{X}}=\sum_{i=1}^r \binom{1+\sum_j n_j}{1+n_i}L_i^{1+n_i}c(n_1,\dots,\widehat{n_i},\dots,n_r)\prod_{j\neq i}(L_j)_t^{n_j}$$
where $t\in T$ is any closed point. Put $d_i(n_1,\dots,n_r):=\binom{1+\sum_j n_j}{1+n_i}c(n_1,\dots,\widehat{n_i},\dots,n_r)$ to conclude.
\end{proof}

\begin{lemma}\label{thm:Numerical lemma}
	Let $G>0, H\geq 0$ and $l\in (0;1)$ be positive real numbers, and $d\geq 2$ be an integer. 
	Then for every choice of real numbers $a_0,\dots,a_d$ satisfying
			$$\max\left\{-\frac{Gl}{4}, -\frac{G^2}{4H(d-1)}\right\}< a_0\leq 0, \quad a_1\geq G, \quad\text{and}\quad |a_i|\leq H \;\forall i\geq 2,$$
	the polynomial $p(t)=\sum_{i=0}^d a_it^i$ takes a strictly positive value in the interval $(0;l/2)$. (If $H=0$ we set $\frac{G^2}{4H(d-1)}=+\infty$).
	\end{lemma}
	\begin{proof}
	Let $a_0,\dots ,a_d$ be real numbers satisfying the prescribed conditions. 
	 We have, for $0<t<1$:
			\begin{eqnarray*}
				p(t) & = & a_0+a_1t+\sum_{i\geq 2}a_i t^i \\
				& \geq & a_0+Gt-H(d-1)t^2.
			\end{eqnarray*}
	So it it enough to prove that $q(t):=a_0+Gt-H't^2$ takes a strictly positive value on $(0,l/2)$, with $H':=H(d-1)$. First consider the special case where $H=0$. Then $a_0>-Gl/4$, so $q(l/3)>Gl/12>0$. From now we assume that $H>0$. We have to show that $q(t)$ has a real positive root $t_0\in (0,l/2)$ such that $q'(t_0)>0$. Real roots exist if
			\begin{equation*}
			a_0>\frac{-G^2}{4H'}.
			\end{equation*}
which holds by assumption on $a_0$. Then the smallest positive root of $q(t)$ is
		$$t_0=\frac{G-\sqrt{G^2+4a_0H'}}{2H'}.$$
Note that 
		$$q'(t_0)=G-2H't_0=\sqrt{G^2+4a_0H'}> 0.$$
Hence we just have to verify that $t_0<l/2$. This condition is equivalent to 
		$$G-lH'<\sqrt{G^2+4a_0H'}.$$
This inequality is trivially satisfied if $G-lH'<0$. If $G-lH'\geq 0$,  then it is equivalent to
		$$\frac{l(lH'-G)}{4}-\frac{Gl}{4}<a_0,$$
which holds because $l(lH'-G)/4<0$ and $-Gl/4<a_0$ by assumption. Therefore $q(t)$ takes a strictly positive value in $(0,l/2)$, as desired.
	\end{proof}

\bibliographystyle{alpha}
\bibliography{Bibliography}

\end{document}